\documentclass[a4paper,draft]{amsart}
 \usepackage{a4}
 
\usepackage{amsmath,mathrsfs,amssymb,amsthm,amscd}
\usepackage[]{fontenc}
\usepackage{xy} 
\usepackage{enumerate}
\usepackage[hyperindex]{hyperref}
\xyoption{all}

\numberwithin{equation}{section}
\newcommand{\nnpar}[1]{\vskip 10pt \noindent {\scshape #1}}\nopagebreak%
\theoremstyle{plain}
\newtheorem{theorem}[equation]{Theorem}
\newtheorem{corollary}[equation]{Corollary}
\newtheorem{proposition}[equation]{Proposition}
\newtheorem{lemma}[equation]{Lemma}

\theoremstyle{definition}
\newtheorem{definition}[equation]{Definition}
\newtheorem{notation}[equation]{Notation}

\newtheorem{remark}[equation]{Remark}

\DeclareMathOperator{\dv}{div}
\DeclareMathOperator{\diva}{\widehat{div}}

\DeclareMathOperator{\cha}{\widehat{CH}}

\DeclareMathOperator{\CH}{CH}

\DeclareMathOperator{\za}{\widehat{Z}}
\DeclareMathOperator{\rata}{\widehat{Rat}}

\DeclareMathOperator{\dd}{d}
\DeclareMathOperator{\Img}{Im}

\DeclareMathOperator{\cl}{cl}

\DeclareMathOperator{\Ker}{Ker}

\DeclareMathOperator{\Id}{id}

\DeclareMathOperator{\Spec}{Spec}

\DeclareMathOperator{\amap}{a}

\DeclareMathOperator{\fflat}{flat}

\newcommand{\smooth}{\text{{\rm smooth}}}

\newcommand{\alg}{\text{{\rm alg}}}

\newcommand{\extd}{\text{{\rm extd}}}

\newcommand{\D}{\text{{\rm cur}}}

\def\?{\ ???\ \immediate\write16{}%
  \immediate\write16{Warning: There was still a question mark . . . }%
  \immediate\write16{}}

\begin{document}










\title{Semipurity of tempered Deligne cohomology}
\author{Jos\'e Ignacio Burgos Gil}
\address{Facultad de Matem\'aticas\\Universidad de Barcelona\\
Gran V\'\i{}a C. C. 585\\ Barcelona, Spain}
\thanks{Partially supported by Grants
   BFM2003-02914 and MTM2006-14234-C02-01}

 \begin{abstract}
   In this paper we define the formal and tempered Deligne cohomology
   groups, that are obtained by applying the Deligne complex functor
   to the complexes of formal differential forms and tempered currents
   respectively. We then prove the existence of a duality between them,
   a vanishing theorem for the former and
   a semipurity property for the latter. The motivation of these
   results comes from the 
   study of covariant arithmetic Chow groups. The semi-purity property
   of tempered Deligne cohomology implies, in particular, that several
   definitions of covariant arithmetic Chow groups agree for
   projective arithmetic varieties. 
 \end{abstract}

\maketitle

\tableofcontents

\section{Introduction}
\label{sec:introduction}

The aim of this note is to study some properties of formal and
tempered Deligne cohomology (with real coefficients). These cohomology
groups are defined by 
applying the Deligne complex functor to the complexes of formal
differential forms and tempered currents respectively.

Let $X$ be a
complex projective manifold and let $W$ be a Zariski locally closed
subset of $X$. Let $i:W\longrightarrow X$ denote the inclusion and let
$i^{*},i^{!}, i_{\ast}, i_{!}$ be the induced functors in the derived
category of abelian sheaves. Then
the complex of formal differential forms of $W$ 
computes the cohomology of $W$ with compact supports. That is, it
computes the groups $H^{\ast}(X,i_{!}i^{\ast}\underline
{\mathbb{R}})$. The 
complex of tempered currents on $W$ compute the cohomology of $X$
with 
supports on $W$, that is, it computes the groups 
$H^{\ast}(X,i_{*}i^{!}\underline {\mathbb{R}})$. Following Deligne, the
previous groups have a mixed Hodge structure, hence a Hodge
filtration that we will call the Deligne-Hodge filtration. The
complexes of formal differential forms and tempered 
currents are examples of Dolbeault complexes (see
\cite{BurgosKramerKuehn:cacg}). Therefore they have a Hodge
filtration obtained from the bigrading of differential forms. In
general, this Hodge filtration does not induce the Deligne-Hodge
filtration in cohomology. Moreover, the spectral sequence associated
to this Hodge filtration does not degenerate at the $E^{1}$-term. 

This implies that formal and tempered Deligne cohomology groups
with real coefficients will not have, in general, the same properties
as Deligne-Beilinson cohomology. For instance they do not need to be finite
dimensional. They have a structure of topological vector spaces,
but they may be non-separated.

Note however that, in the particular case when $W=X$, the formal and
tempered Deligne cohomology groups  with real coefficients, agree with
the usual real Deligne cohomology groups. 

In this note we will construct a (Poincar\'e like) duality between
formal Deligne 
cohomology and tempered Deligne  cohomology, that induce a perfect
pairing between the corresponding separated vector spaces. In particular,
applying this duality to the case $W=X$ we obtain an exceptional duality for
real Deligne Beilinson cohomology (corollary \ref{cor:3}) of smooth
projective varieties that, to my knowledge, is new. The shape of this
exceptional duality reminds very much the functional equation of
$L$-functions. It would be interesting to know whether this duality
has any arithmetic meaning.

The second result is a vanishing result for formal Deligne cohomology.
Thanks to the previous duality, the vanishing result of formal Deligne
cohomology implies a semipurity property of
tempered Deligne cohomology (corollary 
\ref{cor:1}). 

The motivation for these results comes from the study of covariant
arithmetic Chow groups introduced in \cite{Burgos:acr} and
\cite{BurgosKramerKuehn:cacg}. The
covariant arithmetic Chow groups are a variant of the arithmetic Chow
groups defined by Gillet and Soul\'e, that are covariant for 
arbitrary proper morphism. By contrast, the groups defined by
Gillet and Soul\'e are only covariant for proper morphisms
between arithmetic varieties that induce  smooth  maps between the
corresponding complex varieties. The covariant arithmetic Chow groups
do not have a product structure, but they are a module over the
contravariant arithmetic Chow groups (see
\cite{BurgosKramerKuehn:cacg} for more details). Similar definitions
of covariant Chow groups have been given by  Kawaguchi and Moriwaki
\cite{KawaguchiMoriwaki:isfav} 
and by Zha \cite{zha99:_rieman_roch}. These two definitions are
equivalent except for the fact that Zha neglects the structure of real
manifold induced on the complex manifold associated to an arithmetic variety.

Although not explicitly stated, in the paper
\cite{BurgosKramerKuehn:cacg}, the covariant arithmetic 
Chow groups are defined by means of  tempered Deligne
cohomology. The semi-purity 
property of tempered Deligne cohomology was announced and used in
\cite{BurgosKramerKuehn:cacg}. Hence this paper can be seen as a
complement of \cite{BurgosKramerKuehn:cacg}.  A new consequence of the
semipurity 
property is that, for an arithmetic variety that is generically
projective, the covariant Chow groups introduced in \cite{Burgos:acr}
and \cite{BurgosKramerKuehn:cacg}
are isomorphic to the covariant Chow groups introduced
by Kawaguchi and Moriwaki.

\nnpar{Acknowledgments.} In the course of preparing this manuscript,
I had many stimulating discussions with many colleagues. We would 
like to thank them all. In particular, I would like to express my 
gratitude to J.-B.~Bost, U.~K\"uhn, J.~Kramer, K.~K\"uhnemann, V.~Maillot, 
D.~Roessler, C.~Schapira, J~Wildeshaus. Furthermore, I would 
like to thank the CRM (Bellaterra, Barcelona), for partial support of
this work.

\section{Complexes of forms and currents}
\label{sec:compl-forms-curr}

By a complex algebraic manifold we will mean the analytic manifold
associated to a smooth scheme over $\mathbb{C}.$ 
Let $X$ be a projective complex algebraic manifold. We will consider
the following situation:  let $Z\subset Y$ 
be closed subvarieties of $X$, let $U$ and $V$ be the open subsets
$U=X\setminus Y$, $V=X\setminus
Z$ and let $W$ be the locally closed subset $W=Y\setminus Z$.

\subsection{Flat forms and Whitney forms}
\label{sec:flat-forms-whitney}
\ 

\nnpar{The complex of Whitney forms.} Let
$\mathscr{E}^{\ast}_{X}$ denote the 
sheaf of smooth differential forms on $X$. We will denote by $E^{\ast}(U)$ the complex
of global differential forms over $U$  and by $E^{\ast}_{c}(U)$ the complex of
differential forms with compact support.

Let
$\mathscr{E}^{\ast}_{X}(\fflat Y)$ denote the ideal sheaf of
differential forms that are flat
along $Y$. Recall that a differential form on $X$ is called flat
along $Y$ if its
Taylor expansion vanishes at all points of
$Y$. We write
$$\mathscr{E}^{\ast}_{Y^{\infty}}=
\mathscr{E}^{\ast}_{X}/\mathscr{E}^{\ast}_{X}(\fflat 
Y).$$
The sections of this complex of sheaves are called Whitney forms on $Y$.
Whitney's extension theorem (\cite{Tougeron:Ifd} IV theorem 3.1),
gives us a precise description of the space of Whitney 
forms in terms of jets over $Y$. For instance, if $Y$ is the smooth
subvariety of $\mathbb{C}^{n}$ defined by the equations $z_1=\dots
=z_{k}=0$, then the germ of the sheaf of Whitney functions on  $Y$ at
the point $x=(0,\dots ,0)$ is 
\begin{displaymath}
  \mathscr{E}^{0}_{Y^{\infty},x}=
  \mathscr{E}^{0}_{Y,x}[[z_{k+1},\dots ,z_{n},\bar z_{k+1},\dots ,\bar
  z_{n}]]. 
\end{displaymath}

We will 
write 
\begin{displaymath}
  \mathscr{E}^{\ast}_{Y^{\infty}}(\fflat Z)=
\mathscr{E}^{\ast}_{X}(\fflat Z)/\mathscr{E}^{\ast}_{X}(\fflat 
Y).
\end{displaymath}
Observe that $\mathscr{E}^{\ast}_{Y^{\infty}}(\fflat Z)$ can also be
defined as the kernel of the morphism 
\begin{displaymath}
  \mathscr{E}^{\ast}_{Y^{\infty}}\longrightarrow
  \mathscr{E}^{\ast}_{Z^{\infty}}. 
\end{displaymath}

The sheaf $\mathscr{E}^{\ast}_{Y^{\infty}}(\fflat Z)$
agrees with the sheaf denoted $\mathbb{C}_{W}\overset{W}{\otimes }
\mathcal{C}_{X}^{\infty} $ in \cite{KashiwaraSchapira:mfcacs}.

The complex 
$\mathscr{E}^{\ast}_{Y^{\infty}}(\fflat Z)$ is a complex of fine sheaves. 
We will denote the corresponding 
complex of global sections by
$E^{\ast}_{X^{\mathcal{W}}}(W):=\Gamma
(X,\mathscr{E}^{\ast}_{Y^{\infty}}(\fflat Z))$.  Note that the complex
$E^{\ast}_{X^{\mathcal{W}}}(W)$ depends only on the locally closed
subspace $W\subset X$ and not on a particular choice of closed subsets $Y$ and
$Z$. Observe also that
$E^{\ast}_{X^{\mathcal{W}}}(X)=E^{\ast}(X)$ is the usual complex of smooth
differential forms on $X$.

We will denote by
$E^{\ast}_{X^{\mathcal{W}},\mathbb{R}}(W)$ the
real subcomplex underlying  $E^{\ast}_{X^{\mathcal{W}}}(W)$.

By the 
acyclicity of fine sheaves, there is a diagram of short exact
sequences   
\begin{equation}\label{eq:1}
  \xymatrix{ &&&0\ar[d]&\\
    & 0 \ar[d]&  0 \ar[d]& E_{X^{\mathcal{W}}}^{\ast}(W) \ar[d]&\\
    0\ar[r] &E_{X^{\mathcal{W}}}^{\ast}(U)\ar[r]\ar[d]&
    E^{\ast}(X)\ar[r]\ar[d]&
    E_{X^{\mathcal{W}}}^{\ast}(Y)\ar[r]\ar[d]&0\\
    0\ar[r] &E_{X^{\mathcal{W}}}^{\ast}(V)\ar[r]\ar[d]&
    E^{\ast}(X)\ar[r]\ar[d]&
    E^{\ast}_{X^{\mathcal{W}}}(Z)\ar[r]\ar[d]&0\\
    & E_{X^{\mathcal{W}}}^{\ast}(W) \ar[d]& 0 & 0 & \\
    & 0 &&&
}
\end{equation}

The complex $E^{\ast}(X)$ is a topological vector space with the
$C^{\infty}$ topology. With this topology $E^{\ast}(X)$ is a Fr\'echet
topological vector space (\cite{bourbaki87:_topol_vector_spaces_chapt}
III p. 9). Moreover $E^{\ast}_{X^{\mathcal{W}}}(U)$ is a
closed subspace.  In fact, by \cite{Tougeron:Ifd} V corollaire 1.6,
it is the closure of  the 
complex of differential forms that have compact support contained in
$U$, that we denote $E^{\ast}_{c}(U)$. More generally, all the
monomorphisms in 
diagram \eqref{eq:1} are closed immersions.

The following result states that, being $U$ an algebraic open subset of $X$,
the complex $E^{\ast}_{X^{\mathcal{W}}}(U)$ 
does not depend on $X$ but only on $U$. 

\begin{proposition}\label{prop:3}
Let $\pi :\widetilde X\longrightarrow X$ be a proper birational
morphism with $D=\pi ^{-1}(Y)$, that induces an isomorphism between
$\widetilde X\setminus D$ and $U$, then the natural map  
\begin{displaymath}
  \pi ^{\ast}:E^{\ast}(X)
  \longrightarrow 
  E^{\ast}(\widetilde X)
\end{displaymath}
induces an isomorphism $\pi ^{\ast}:\Gamma (X,\mathscr{E}^{\ast}_{X}(\fflat Y))
  \longrightarrow 
  \Gamma (\widetilde X,\mathscr{E}^{\ast}_{\widetilde X}(\fflat D))$. 
\end{proposition}
\begin{proof}
  By \cite{Poly:shcsesa}  the morphism
  \begin{displaymath}
    \pi ^{\ast}:E^{\ast}(X)\longrightarrow E^{\ast}(\widetilde X)
  \end{displaymath}
  is a closed immersion. Since $\Gamma
  (X,\mathscr{E}^{\ast}_{X}(\fflat Y)) 
  $ and $
  \Gamma (\widetilde X,\mathscr{E}^{\ast}_{\widetilde X}(\fflat D))$
  are the closure of $E^{\ast}_{c}(U)$ in $E^{\ast}(X)$ and
  $E^{\ast}(\widetilde X)$ respectively, then they are identified by
  $\pi ^{\ast}$.  
\end{proof}

\nnpar{The cohomology of the complex of Whitney forms.} By
\cite{Poly:shcsesa} (see also \cite{BrasseletPflaum:_whitn} for a
more general statement) we have

\begin{proposition} \label{prop:4}
  The complex $\mathscr{E}^{\ast}_{Y^{\infty}}$ is a resolution
  of the constant 
  sheaf $\underline {\mathbb{C}}$ on $Y$ by fine sheaves. Therefore
  \begin{displaymath}
    H^{\ast}(E_{X^{\mathcal{W}}}^{\ast}(W))=H^{\ast}_{c}(W,\mathbb{C}),
  \end{displaymath}
  where $H^{\ast}_{c}$ denotes cohomology with compact supports.
  \hfill $\square$
\end{proposition}

\subsection{Currents with support in a subvariety}
\label{sec:curr-with-supp}
\ 

\nnpar{The complex of currents.} We first recall the definition of the
complex of currents and we fix the sign convention and some
normalizations. We 
will follow the conventions of \cite{BurgosKramerKuehn:cacg} \S 5.4
but with the homological grading. 
 
 Let $\mathscr
{D}_{n}^{X}$ be the sheaf of degree $n$ currents on $X$. That is,
for any open subset $V$ of $X$, the group $\mathscr 
{D}_{n}^{X}(V)$ is the topological dual of the group of sections 
with compact support $E^{n}_{c}(V)$. The 
differential 
\begin{displaymath}
\dd:\mathscr{D}_{n}^{X}\longrightarrow\mathscr{D}_{n-1}^{X}
\end{displaymath}
is defined by
\begin{displaymath}
\dd T(\varphi)=(-1)^{n}T(\dd\varphi);
\end{displaymath}
here $T$ is a current and $\varphi$ a test form.
Note that we are using the sign convention of, for instance
\cite{Jannsen:DcHD},  instead of the sign convention of
\cite{GriffithsHarris:pag}. 

The bigrading
$\mathscr{E}^{n}_{X}=\bigoplus_{p+q=n}\mathscr{E}^{p,q}_{X}$ induces a
bigrading
$$\mathscr{D}_{n}^{X}=\bigoplus_{p+q=n}\mathscr{D}_{p,q}^{X},$$ with 
$\mathscr 
{D}_{p,q}^{X}(V)$ the topological dual of
$\Gamma_{c}(V,\mathscr{E}^{p,q}_{X})$. 

The real structure of $\mathscr{E}^{n}_{X}$ induces a real structure 
\begin{displaymath}
  \mathscr{D}_{n}^{X,\mathbb{R}}\subset \mathscr{D}_{n}^{X}.
\end{displaymath}

We will denote 
\begin{displaymath}
  \mathscr{D}_{n}^{X,\mathbb{R}}(p)=\frac{1}{(2\pi i)^{p}}
  \mathscr{D}_{n}^{X,\mathbb{R}}\subset \mathscr{D}_{n}^{X}. 
\end{displaymath}

If $X$ is equidimensional of dimension $d$ we will write
\begin{equation}\label{eq:4}
  \mathscr{D}^{n}_{X}=\mathscr{D}_{2d-n}^{X}, \quad
  \mathscr{D}^{p,q}_{X}=\mathscr{D}_{d-p,d-q}^{X}, \quad
  \text{and}\quad
  \mathscr{D}^{n}_{X,\mathbb{R}}(p)=\mathscr{D}_{2d-n}^{X,\mathbb{R}}(d-p). 
\end{equation}

We will use all the conventions of \cite{BurgosKramerKuehn:cacg} \S
5.4. In particular, if $y$ is an algebraic cycle of $X$ of dimension
$e$, we will write $\delta _{y}\in \mathscr{D}_{e,e}^{X}\cap
\mathscr{D}_{2e}^{X,\mathbb{R}}(e)$ for the current 
\begin{displaymath}
  \delta _{y}(\eta)=\frac{1}{(2\pi i)^{e}}\int_{y}\eta.  
\end{displaymath}
Furthermore, there is an action
\begin{displaymath}
  \begin{matrix}
    \mathscr{E}^{n}_{X}\otimes\mathscr{D}_{m}^{X}&\longrightarrow &
    \mathscr{D}_{m-n}^{X},\\ 
    \omega\otimes T&\longmapsto &\omega 
    \land T    
  \end{matrix}
\end{displaymath}
where the current $\omega \land T$ is defined by
\begin{displaymath}
(\omega\land T)(\eta)=T(\eta\land\omega).
\end{displaymath}
This action induces actions
\begin{displaymath}
  \mathscr{E}^{p,q}_{X}\otimes\mathscr{D}_{r,s}^{X}\longrightarrow 
\mathscr{D}_{r-p,s-q}^{X}, \quad\text{and}\quad 
\mathscr{E}^{n}_{X,\mathbb{R}}(p)\otimes
\mathscr{D}_{m}^{X,\mathbb{R}}(q)\longrightarrow   
\mathscr{D}_{m-n}^{X,\mathbb{R}}(q-p).
\end{displaymath}

Finally, if $X$ is equidimensional of dimension $d$, there is a
fundamental current $\delta _{X}\in \mathscr{D}_{d,d}^{X}\cap
\mathscr{D}_{2d}^{X,\mathbb{R}}(d)$, and a morphism 
\begin{equation}\label{eq:11}
\mathscr{E}^{\ast}_{X}\longrightarrow \mathscr{D}_{2d-\ast}^{X}= 
\mathscr{D}^{\ast}_{X},\quad\omega \longmapsto [\omega ]=\omega 
\land \delta _{X}.
\end{equation}
This morphism sends $\mathscr{E}^{n}_{X\mathbb{R}}(p)$ to
$\mathscr{D}_{2d-n}^{X,\mathbb{R}}(d-p)=\mathscr{D}^{n}_{X,\mathbb{R}}(p)$.

\nnpar{Currents with support on a subvariety and tempered currents.} 
As in the previous section let $Z\subset Y$ denote two closed
subvarieties of $X$ 
and put $U=X\setminus Y$, $V=X\setminus Z$ and $W=Y\setminus Z$.  
We denote by $\mathscr
{D}_{\ast}^{Y^{\infty}}$ the subcomplex of $\mathscr{D}_{\ast}^{X}$
formed by currents with support on
$Y$. In other words, for any open subset $U'$ of $X$ we have
\begin{displaymath}
  \mathscr
{D}_{n}^{Y^{\infty}}(U')=\{T\in \mathscr
{D}_{n}^{X}(U')\mid T(\eta)=0,\ \forall \eta\in
\Gamma_{c}(U'\cap U,\mathscr{E}^{n}_{X})\}. 
\end{displaymath}
Observe that, by continuity, the sections of
$\mathscr{D}_{n}^{Y^{\infty}}(U')$ vanish on the subgroup $\Gamma
_{c}(U',\mathscr{E}^{\ast}_{X}(\fflat Y))$. 

We write $\mathscr{D}_{n}^{X/Y^{\infty}}=\mathscr
{D}_{n}^{X}\left / \mathscr
{D}_{n}^{Y^{\infty}}\right .$ and $
  \mathscr{D}_{n}^{Y^{\infty}/Z^{\infty}}=
  \mathscr{D}_{n}^{Y^{\infty}}/\mathscr{D}_{n}^{Z^{\infty}}.$

As in the case of differential forms, the complex $
\mathscr{D}_{n}^{Y^{\infty}/Z^{\infty}}$ can also be defined as the
kernel of the morphism 
\begin{displaymath}
   \mathscr{D}_{n}^{X/ Z^{\infty}} \longrightarrow
   \mathscr{D}_{n}^{X/ Y^{\infty}}. 
\end{displaymath}
All the above sheaves inherit a bigrading and a real structure.

Observe that, except for the fact that we are using here the
homological grading, the complex of sheaves $\mathscr{D}_{n}^{X/ Y^{\infty}}$
agrees with the complex denoted by
$\mathcal{TH}om(\mathbb{C}_{W},\mathcal{D}b_{X})$ in
\cite{KashiwaraSchapira:mfcacs}.

The complex $\mathscr{D}_{n}^{ Y^{\infty}/Z^{\infty}}$ is a complex of fine
sheaves. We will denote the complex of global sections by
$D_{\ast}^{X^{\mathcal{T}}}(W^{\infty})=\Gamma (X,
\mathscr{D}_{\ast}^{Y^{\infty}/Z^{\infty}})$. Thus the complex
$D_{\ast}^{X^{\mathcal{T}}}(W^{\infty})$ is defined for any Zariski
locally closed subset $W\subset X$. The corresponding real
complex will be denoted by $D_{\ast}^{X^{\mathcal{T}},\mathbb{R}}(W^{\infty})$.

By \cite{Poly:shcsesa}, the complex $D_{\ast}^{X^{\mathcal{T}}}(U)$ can be
identified with the
image of the morphism
\begin{displaymath}
  D^{\ast}(X)\longrightarrow D^{\ast}(U).
\end{displaymath}
That is, it is the complex of currents on $U$ that can be extended to a
current on the whole $X$. The elements of $D_{\ast}^{X^{\mathcal{T}}}(U)$ will
be called tempered currents. In the literature they are called also
moderate, temperate or extendable currents. Moreover, as was the case with the complex
$E^{\ast}_{X^{\mathcal{W}}}(U)$, being $U$ a Zariski open subset, the
complex $D_{\ast}^{X^{\mathcal{T}}}(U)$ 
only depends on $U$ 
and not on $X$. 

\nnpar{The pairing between forms and currents.} We have already
introduced an action  
\begin{equation}\label{eq:5}
E^{n}(X)\otimes D_{m}(X)\longrightarrow 
D_{m-n}(X),\quad\omega\otimes T\longmapsto\omega 
\land T,
\end{equation}
where the current $\omega \land T$ is defined by
\begin{displaymath}
(\omega\land T)(\eta)=T(\eta\land\omega).
\end{displaymath}

The subspace $D_{\ast}^{X^{\mathcal{T}}}(Y)$ is invariant under this action and
annihilates 
the subspace $E_{X^{\mathcal{W}}}^{\ast}(U)$. Therefore we obtain induced actions
\begin{equation}\label{eq:6}
  E^{n}_{X^{\mathcal{W}}}(Y)\otimes D^{X^{\mathcal{T}}}_{m}(Y)\longrightarrow 
D^{X^{\mathcal{T}}}_{m-n}(Y),\qquad E^{n}_{X^{\mathcal{W}}}(U)\otimes
D_{m}^{X^{\mathcal{T}}}(U)\longrightarrow  
D_{m-n}^{X^{\mathcal{T}}}(U)
\end{equation}
and, more generally, an action
\begin{equation}\label{eq:7}
    E^{n}_{X^{\mathcal{W}}}(W)\otimes
    D_{m}^{X^{\mathcal{T}}}(W)\longrightarrow  
    D_{m-n}^{X^{\mathcal{T}}}(W). 
\end{equation}

 Since $X$ is proper, there is a canonical morphism
\begin{displaymath}
  \deg: D_{0}(X)\longrightarrow \mathbb{C}
\end{displaymath}
given by $\deg(T)=T(1)$. Observe that
$\deg(D_{0}^{\mathbb{R}}(X))\subset \mathbb{R}$.

Combining the degree and the above actions, we recover the pairing
\begin{displaymath}
  E^{n}(X)\otimes D_{n}(X)\longrightarrow \mathbb{C},
\end{displaymath}
that identifies $D_{n}(X)$ with the topological dual of
$E^{n}(X)$. Under this identification, the subspace $E^{n}_{X^{\mathcal{W}}}(U)$
is the orthogonal to the subspace $D_{n}^{X^{\mathcal{T}}}(Y)$. Therefore
$D_{n}^{X^{\mathcal{T}}}(U)$ is the topological dual of
$E^{n}_{X^{\mathcal{W}}}(U)$  
and $D_{n}^{X^{\mathcal{T}}}(Y)$ is the topological dual of
$E^{n}_{X^{\mathcal{W}}}(Y)$. More generally $D_{n}^{X^{\mathcal{T}}}(W)$ is
the topological dual of $E^{n}_{X^{\mathcal{W}}}(W)$. Note that
here, the key point is the fact that $E^{n}_{X^{\mathcal{W}}}(U)$ is the closure
of $\Gamma _{c}(U,\mathscr{E}^{n}_{X})$ and hence a closed subspace. 

The above pairings induce a pairing
\begin{displaymath}
    E^{n}_{\mathbb{R}}(X)(p)\otimes
    D_{n}^{\mathbb{R}}(X)(p)\longrightarrow \mathbb{R}, 
\end{displaymath}
and similar pairings for the other complexes of forms and currents.

Finally, observe that there is a commutative diagram with exact rows
and columns 
\begin{equation}\label{eq:2}
  \xymatrix{ &&&0\ar[d]&\\
    & 0 \ar[d]&  0 \ar[d]& D^{X^{\mathcal{T}}}_{\ast}(W) \ar[d]&\\
    0\ar[r] &D_{\ast}^{X^{\mathcal{T}}}(Z)\ar[r]\ar[d]&
    D_{\ast}(X)\ar[r]\ar[d]&
    D_{\ast}^{X^{\mathcal{T}}}(V)\ar[r]\ar[d]&0\\
    0\ar[r] &D_{\ast}^{X^{\mathcal{T}}}(Y)\ar[r]\ar[d]&
    D_{\ast}(X)\ar[r]\ar[d]&
    D_{\ast}^{X^{\mathcal{T}}}(U)\ar[r]\ar[d]&0\\
    & D^{X^{\mathcal{T}}}_{\ast}(W) \ar[d]& 0 & 0 & \\
    & 0 &&&
}
\end{equation}
that is the topological dual of the diagram \eqref{eq:1}.

\nnpar{The homology of the complexes of currents.} By \cite{Poly:shcsesa}
we have 
\begin{proposition}
  The homology of the complexes $D_{\ast}^{X^{\mathcal{T}}}(W)$ is
  given by 
  \begin{displaymath}
     H_{\ast}(D^{X^{\mathcal{T}}}_{\ast}(W))=H_{\ast}^{BM}(W,\mathbb{C}),
  \end{displaymath}
  where $H_\ast^{BM}$ denote Borel-Moore homology. In particular,
  since we are assuming $Y$ proper,
  \begin{displaymath}
    H_{\ast}(D_{\ast}^{X^{\mathcal{T}}}(Y))=H_{\ast}(Y,\mathbb{C}).
  \end{displaymath}
\hfill $\square$
\end{proposition}

\subsection{Formal and tempered Deligne cohomology}
\label{sec:form-deligne-cohom}
\

\nnpar{Formal Deligne cohomology.} 
The complex 
$E^{\ast}_{X^{\mathcal{W}},\mathbb{R}}(W)$ is an example of a Dolbeault
algebra (see \cite{BurgosKramerKuehn:cacg}). 
Recall that, following Deligne, the cohomology of any
complex variety has a mixed 
Hodge structure. We will call the Hodge filtration of this mixed
Hodge structure the Deligne-Hodge filtration.

From the structure of Dolbeault algebra of
$E^{\ast}_{X^{\mathcal{W}}}(W)$ we can define a Hodge
filtration. It is the filtration associated to the bigrading. 
In general, this Hodge filtration does not induce the
Deligne-Hodge filtration in cohomology. Moreover, the spectral sequence
associated to this Hodge filtration does not need to degenerate at the
$E_{1}$ term. Therefore, the Dolbeault cohomology groups
$H^{p,q}_{\overline \partial}(E^{\ast}(Y^{\infty}))$ are not, in general, direct
summands of $H^{p+q}(Y,\mathbb{C})$. In fact, they can be infinite
dimensional as can be seen in the easiest example: Put
$X=\mathbb{P}^{1}_{\mathbb{C}}$. Let $t$ be the absolute coordinate
and let $Y$ be
the point $t=0$. Then $H^{0,0}_{\bar
  \partial}(E^{\ast}(Y^{\infty}))=\mathbb{C}[[t]]$, the ring of formal
power series in one variable.

Following \cite{Burgos:CDB} and \cite{BurgosKramerKuehn:cacg}, to
every Dolbeault algebra we can associate a Deligne algebra. 
We
refer the reader to \cite{Burgos:CDB} and
\cite{BurgosKramerKuehn:cacg} \S 5 for the definition and properties
of Dolbeault algebras, Dolbeault complexes and the associated Deligne
complexes. We will use freely the notation therein. In particular the
Deligne algebra associated to the above Dolbeault algebra
will be denoted $\mathcal{D}^{\ast}(E^{\ast}_{X^{\mathcal{W}}}(W),\ast)$.

\begin{definition} \label{def:1} The real formal Deligne cohomology
  of $W$ (with
  compact supports) is defined by
  \begin{align*}
    H_{\mathcal{D}^{f},c}^{\ast}(W^{\infty},\mathbb{R}(p))&=
    H^{\ast}(\mathcal{D}^{\ast}(E_{X^{\mathcal{W}}}(W),p)).
  \end{align*}
  When $W$ is proper we will just write
  $H_{\mathcal{D}^{f}}^{\ast}(W^{\infty},\mathbb{R}(p))$. 
\end{definition}
The notation $W^{\infty}$ is a reminder that this cohomology depends,
not only on $W$ but on an infinitesimal neighborhood of infinite
order of $W$ in $X$.  

\begin{remark}
  Since we are assuming that $X$ is smooth and proper, the formal Deligne
  cohomology of $X$,
  $H_{\mathcal{D}^{f}}^{\ast}(X^{\infty},\mathbb{R}(p))$, given in
  the previous definition, agrees with the usual Deligne cohomology of $X$.
  Nevertheless, by the discussion before the definition, the formal
  Deligne cohomology with 
  compact supports of $U$ or the formal Deligne cohomology of $Y$,
  do not agree, in general, with the usual Deligne-Beilinson
  cohomology. For instance the groups
  $H_{\mathcal{D}^{f}}^{\ast}(U,\mathbb{R}(p))$ can be infinite
  dimensional. 
\end{remark}

\nnpar{Homological Dolbeault complexes and homological Deligne
  complexes.} In order to define formal Deligne homology we 
first translate the notions of \cite{BurgosKramerKuehn:cacg}
\S 5.2 to the homological grading.

\begin{definition}
\label{def:12}
A \emph{homological Dolbeault complex} $A=(A_{\ast}^{\mathbb{R}},\dd_{A})$ is 
a graded complex of real vector spaces, which is bounded from above 
and equipped with a bigrading on $A^{\mathbb{C}}=A^{\mathbb{R}}
\otimes_{\mathbb{R}}{\mathbb{C}}$, i.e.,
\begin{displaymath}
A_{n}^{\mathbb{C}}=\bigoplus_{p+q=n}A_{p,q},
\end{displaymath}  
satisfying the following properties:
\begin{enumerate}
\item[(i)]
The differential $\dd_{A}$ can be decomposed as the sum $\dd_{A}=
\partial+\bar{\partial}$ of operators $\partial$ of type $(-1,0)$, 
resp. $\bar{\partial}$ of type $(0,-1)$.
\item[(ii)] 
It satisfies the symmetry property $\overline{A_{p,q}}=A_{q,p}$,
where $\overline{\phantom{M}}$ denotes complex conjugation.
\end{enumerate}
\end{definition}

\begin{notation}
\label{def:13}
Given a homological Dolbeault complex $A=(A_{\ast}^{\mathbb{R}},\dd_{A})$, we 
will use the following notations. The Hodge filtration $F$ of $A$ 
is the increasing filtration of $A^{\mathbb{C}}_{\ast}$ given by
\begin{displaymath}
F_{p}A_{n}=F_{p}A_{n}^{\mathbb{C}}=\bigoplus_{p'\leq p}A_{p',n-p'}.
\end{displaymath}
The filtration $\overline F$ of $A$ is the complex conjugate of $F$, 
i.e.,
\begin{displaymath}
\overline{F}_{p}A_{n}=\overline{F}_{p}A_{n}^{\mathbb{C}}=\overline
{F_{p}A_{n}^{\mathbb{C}}}.
\end{displaymath}
For an element $x\in A^{\mathbb{C}}$, we write $x_{i,j}$ for its 
component in $A_{i,j}$. For $k,k' \in \mathbb{Z}$, we define an operator 
$F_{k,k'}:A^{\mathbb{C}}\longrightarrow A^{\mathbb{C}}$ by the 
rule 
\begin{displaymath}
F_{k,k'}(x):=\sum_{l\leq k,l'\leq k'}x_{l,l'}.
\end{displaymath}
We note that the operator $F_{k,k'}$ is the projection of $A^{\ast}_
{\mathbb{C}}$ onto the subspace $F_{k}A_{\ast}\cap\overline{F}_{k'}
A_{\ast}$. This subspace will be denoted $F_{k,k'}A_{\ast}$. We will
also denote by $F_{k}$ the operator $F_{k,\infty}$.  

We denote by $A_{n}^{\mathbb{R}}(p)$ the subgroup $(2\pi i)^{-p}\cdot
A_{n}^{\mathbb{R}}\subseteq A_{n}^{\mathbb{C}}$, and we define the 
operator
\begin{displaymath}
\pi_{p}:A^{\mathbb{C}}\longrightarrow A^{\mathbb{R}}(p)
\end{displaymath}
by setting $\pi_{p}(x):=\frac{1}{2}(x+(-1)^{p}\bar{x})$.
\end{notation}

To any homological Dolbeault complex we can associate a homological
Deligne complex. 

\begin{definition}
Let $A$ be a homological Dolbeault complex. We denote by $A_{\ast}(p)^
{\mathcal{D}}$ the complex $s(A^{\mathbb{R}}(p)\oplus F_{p}A
\overset{u}{\longrightarrow}A^{\mathbb{C}})$, where $u(a,f)=
-a+f$ and $s(\ )$ denotes the simple complex of a morphism of
complexes. 
\end{definition}

\begin{definition}
Let $A$ be a homological Dolbeault complex. Then, the
\emph{(homological) Deligne complex  
$(\mathcal{D}^{\ast}(A,\ast),\dd_{\mathcal{D}})$ associated to $A$} 
is the graded complex given by   
\begin{align*}
&\mathcal{D}_{n}(A,p)=
\begin{cases}
A^{\mathbb{R}}_{n+1}(p+1)\cap F_{n-p,n-p}A_{n+1}^{\mathbb{C}},
&\qquad\text{if}\quad n\geq 2e+1, \\
A^{\mathbb{R}}_n(p)\cap F_{p,p}A_{n}^{\mathbb{C}},
&\qquad\text{if}\quad n\leq 2p,  
\end{cases}  
\intertext{with differential given, for $x\in\mathcal{D}_{n}(A,p)$,  by}
&\dd_{\mathcal{D}}x=
\begin{cases}
-F_{n-p+1,n-p+1}\dd_{A}x,
&\qquad\text{if}\quad n>2p+1, \\
-2\partial\bar{\partial}x, 
&\qquad\text{if}\quad n=2p+1, \\
\dd_{A}x, 
&\qquad\text{if}\quad n\leq 2p.
\end{cases}  
\end{align*}  
\end{definition}

For instance, let $A$ be a Dolbeault complex satisfying $A_{p,q}=0$ 
for $p<0$, $q<0$, $p>n$, or $q>n$.  Then, for $p\ge n$, the complex
$\mathcal{D}(A,p)$  
agrees with the real complex $A_{\ast}^{\mathbb{R}}(p)$. For $0\le p<n$, we
have represented $\mathcal{D}(A,p)$ in figure \ref{fig:1}, where 
the upper right square is shifted by one; this means in particular that
$A_{n,n}$ sits in degree $2n-1$ and $A_{p+1,p+1}$ sits in degree
$2p+1$. For $p<0$ the complex $\mathcal{D}(A,p)$ 
agrees with the real complex $A_{\ast}^{\mathbb{R}}(p+1)[1]$.

\begin{figure}[htb]
\begin{center}
\begin{displaymath}
\xymatrix{&
\makebox{$\begin{pmatrix}
A_{p+1,n}&\leftarrow &\cdots&\leftarrow &A_{n,n} \\
\downarrow &&&& \downarrow \\
\vdots &&&& \vdots \\
\downarrow &&&& \downarrow \\
A_{p+1,p+1}&\leftarrow&\cdots&\leftarrow&A_{n,p+1} 
\end{pmatrix}$}_{\makebox[0pt]{$\mathbb{R}$}}
\makebox[0pt][l]{$(p+1)$} 
\ar[dl]^{-2\partial\overline{\partial}}\\
\makebox{$\begin{pmatrix}
A_{0,p}&\leftarrow&\cdots&\leftarrow&A_{p,p} \\
\downarrow &&&& \downarrow \\
\vdots &&&& \vdots \\
\downarrow &&&& \downarrow \\
A_{0,0}&\leftarrow&\cdots&\leftarrow&A_{p,0}
\end{pmatrix}$}_{\makebox[0pt]{$\mathbb{R}$}}
\makebox[0pt][l]{$(p)$}
&}
\end{displaymath}
\caption{$\mathcal{D}(A,p)$}
\label{fig:1}
\end{center}
\end{figure}

\begin{remark} 
It is clear from the definition that, for all $p\in \mathbb{Z}$, the functor $\mathcal{D}
(\cdot,p)$ is exact.
\end{remark}

The main property of the Deligne complex is expressed by the 
following proposition; for a proof in the cohomological case see
\cite{Burgos:CDB}. 

\begin{proposition}
\label{prop:32}
The complexes $A_{\ast}(p)^{\mathcal{D}}$ and $\mathcal{D}_{\ast}
(A,p)$ are homotopically equivalent. The homotopy equivalences
$\psi:A_{n}(p)^{\mathcal{D}}\longrightarrow\mathcal{D}_{n}(A,p)$,
and $\varphi:\mathcal{D}_{n}(A,p)\longrightarrow A_{n}(p)^{\mathcal 
{D}}$ are given by
\begin{displaymath}
\psi(a,f,\omega)=
\begin{cases}
\pi(\omega),\qquad&\text{if }n\ge 2p+1, \\
F_{p,p}a+2\pi_{p}(\partial\omega_{p+1,n-p-1}),\quad&\text{if }n\le 2p,
\end{cases}
\end{displaymath}  
where $\pi(\omega)=\pi_{p+1}(F_{n-p,n-p}\omega)$, i.e., $\pi$ is
the projection of $A_{\mathbb{C}}$ over the co\-kernel of $u$, and
\begin{displaymath}
\varphi(x)=
\begin{cases}
(\partial x_{p+1,n-p}-\bar{\partial}x_{n-p,p+1},2\partial
x_{p+1,n-p},x),\quad&\text{if }n\ge 2p+1, \\
(x,x,0),&\text{if }n\le 2p.
\end{cases}
\end{displaymath}
Moreover, $\psi\circ\varphi=\Id$, and $\varphi\circ\psi-\Id=\dd h+
h\dd$, where $h:A_{n}(p)^{\mathcal{D}}\longrightarrow A_{n+1}(p)^
{\mathcal{D}}$ is given by
\begin{displaymath}
h(a,f,\omega)=
\begin{cases}
(\pi_{p}(\overline{F}_{p}\omega+\overline{F}_{n-p}\omega),-2F_{p}
(\pi_{p+1}\omega),0),\quad&\text{if }n\ge 2p+1, \\
(2\pi_{p}(\overline{F}_{n-p}\omega),-F_{p,p}\omega-2F_{n-p}
(\pi_{p+1}\omega),0),\quad&\text{if }n\le 2p.
\end{cases}
\end{displaymath}
\end{proposition}
\hfill $\square$

\nnpar{Tempered Deligne homology.} Applying the above discussion to the
complex of currents
$D_{\ast}^{X^{\mathcal{T}},\mathbb{R}}(W)$ we define the homological
Deligne complex 
$\mathcal{D}_{\ast}(D_{\ast}^{X^{\mathcal{T}}}(W),\ast).$

\begin{definition} \label{def:2} \emph{The tempered  Deligne (Borel-Moore)
  homology of $W$} is defined by
  \begin{displaymath}
    H^{\mathcal{D}^{\mathcal{T}}}_{\ast}(W^{\infty},\mathbb{R}(p))=
    H_{\ast}(\mathcal{D}_{\ast}(D^{X^{\mathcal{T}}}_{\ast}(W),p)).
  \end{displaymath}
\end{definition}

\begin{remark}
  \begin{enumerate}
  \item   Again, since $X$ is smooth and proper, the tempered Deligne
    homology of 
    $X$ agrees with the Deligne homology of $X$. In particular,
    the group 
    $H^{\mathcal{D}}_{n}(X,\mathbb{R}(p))$ agrees with the group denoted
    ${}'H^{-n}_{\mathcal{D}}(X,\mathbb{R}(-p))$ in
    \cite{Jannsen:DcHD}. But, since the Hodge filtration of the complex
    of currents with support on $Y$ does not induce the Deligne-Hodge
    filtration in the
    homology of $Y$, the tempered Deligne homology does not
    agree in general with  Deligne-Beilinson homology.
  \item As in the case of formal cohomology, the notation
    $H^{\mathcal{D}^{\mathcal{T}}}_{\ast}(W^{\infty},\mathbb{R}(p))$ reminds us
    that these groups do not depend only on $W$ but on an
    infinitesimal neighborhood of $W$ of infinite order.
  \end{enumerate}
\end{remark}

\nnpar{Equidimensional manifolds.} If $X$ is equidimensional of
dimension $d$ the morphism \eqref{eq:11} induces morphisms
\begin{equation}
    \mathcal{D}^{n}(E^{\ast}(X),p)\longrightarrow 
    \mathcal{D}_{2d-n}(D_{\ast}(X),d-p), \quad p  \in \mathbb{Z},
\end{equation}
that, in turn, induce the Poincar\'e duality isomorphisms
\begin{equation}\label{eq:8}
  H^{n}_{\mathcal{D}}(X,\mathbb{R}(p))
  \longrightarrow H_{2d-n}^{\mathcal{D}}(X,\mathbb{R}(d-p)), \quad n,p \in \mathbb{Z}.
\end{equation}

By analogy, we can define tempered Deligne cohomology groups as follows
\begin{align*}
    H_{\mathcal{D}^{\mathcal{T}}}^{n}(U,\mathbb{R}(p))&=
    H^{\mathcal{D}^{\mathcal{T}}}_{2d-n}(U,\mathbb{R}(d-p)),\\
    H_{\mathcal{D}^{\mathcal{T}},W}^{n}(V,\mathbb{R}(p))&=
    H^{\mathcal{D}^{\mathcal{T}}}_{2d-n}(W^{\infty},\mathbb{R}(d-p)).
\end{align*}

In general, if $X$ is a disjoint union of equidimensional algebraic
manifolds, then we define the tempered Deligne cohomology of $X$
as the direct sum of the tempered Deligne cohomology of its components.

\nnpar{The module structure of tempered Deligne homology.}
The notion of Dolbeault module over a
Dolbeault algebra introduced in \cite{BurgosKramerKuehn:cacg} can be
easily modified to define homological Dolbeault modules over a
Dolbeault algebra. The actions \eqref{eq:5}, \eqref{eq:6} and
\eqref{eq:7} provide the basic examples. Modifying the construction of
\cite{BurgosKramerKuehn:cacg} 5.17 and 5.18 we obtain

\begin{proposition}
  There is a
pseudo-associative action
\begin{displaymath}
    \mathcal{D}^{n}(E_{X^{\mathcal{W}}}(W),p)\otimes 
    \mathcal{D}_{m}(D_{\ast}^{X^{\mathcal{T}}}(W),q)\longrightarrow 
    \mathcal{D}_{m-n}(D_{\ast}^{X^{\mathcal{T}}}(W),q-p)
\end{displaymath}
that induces an associative action
\begin{displaymath}
  H_{\mathcal{D}^{f},c}^{n}(W^{\infty},\mathbb{R}(p))\otimes 
  H^{\mathcal{D}^{\mathcal{T}}}_{m}(W^{\infty},\mathbb{R}(q))\longrightarrow
  H^{\mathcal{D}^{\mathcal{T}}}_{m-n}(W^{\infty},\mathbb{R}(q-p)). 
\end{displaymath}
\hfill $\square$
\end{proposition}

\nnpar{The exceptional duality.} In general, Poincar\'e duality for
Deligne cohomology is not given by a bilinear pairing, but by the
isomorphism \eqref{eq:8} between Deligne cohomology and Deligne
homology (see for instance \cite{Jannsen:DcHD}). Nevertheless, in the case 
of real Deligne cohomology, there is an exceptional duality that comes
from the symmetry of the Deligne complex associated with a Dolbeault
complex. This duality can be generalized to a pairing between formal
Deligne cohomology and tempered Deligne homology.

\begin{proposition}\label{prop:5}
  For every pair of integers $n,p$, there is a pairing
  \begin{displaymath}
    \mathcal{D}^{n}(E_{X^{\mathcal{W}}}(W),p)\otimes
    \mathcal{D}_{n-1}(D^{X^{\mathcal{T}}}(W),p-1)\longrightarrow \mathbb{R}
  \end{displaymath}
  given by $\omega \otimes T\longmapsto T(\omega )$.
  This pairing identifies $\mathcal{D}_{n-1}(D^{X^{\mathcal{T}}}(W),p-1)$
  with the topological dual of
  $\mathcal{D}^{n}(E_{X_{\mathcal{W}}}(W),p)$. Moreover, it is
  compatible, up to the sign, with the differential in the Deligne complex:
  \begin{displaymath}
    T (\dd_{\mathcal{D}}\omega ) =
    \begin{cases}
      (-1)^{n+1}
      (\dd_{\mathcal{D}}T)(\omega),& \text{ if } n\le 2p-1,\\ 
      (-1)^{n}
      (\dd_{\mathcal{D}}T)(\omega),& \text{ if } n\ge 2p.\\ 
    \end{cases}
  \end{displaymath}
  It is also compatible, up to the sign, with the action of
  $\mathcal{D}^{\ast}(E_{X^{\mathcal{W}}}(W),\ast)$. That is, if the
  forms 
  $\omega \in \mathcal{D}^{n}(E_{X^{\mathcal{W}}}(W^{\infty}),p)$ and 
  $\eta \in \mathcal{D}^{l}(E_{X^{\mathcal{W}}}(W),r)$, and the current 
  $T\in \mathcal{D}_{m}(D^{X^{\mathcal{T}}}(W),q)$, with
  $n-m+l=1$ and $p-q+r=1$ then
  \begin{displaymath}
    (\omega \bullet T)(\eta)= 
    \begin{cases}
      (-1)^{n} T(\eta \bullet \omega), & \text{ if } m>2q,\ l\ge 2r,\\       
      T(\eta \bullet \omega), & \text{ if } m\le 2q,\ l< 2r,\\       
      (-1)^{m-1} T(\eta \bullet \omega), & \text{ if } m>2q,\ l< 2r,\\       
      (-1)^{l} T(\eta \bullet \omega), & \text{ if } m\le 2q,\ l\ge 2r.\\       
    \end{cases}
  \end{displaymath}
\end{proposition}
\begin{proof}
  Assume that $n<2p$. Put $q=p-1$ and
  $m=n-1$. Then
  \begin{align*}
    &\mathcal{D}^{n}(E_{X^{\mathcal{W}}}(W),p)\\
    &\phantom{A}=
    E^{n-1}_{X^{\mathcal{W}},\mathbb{R}}(W)(p-1)\left /   
      (F^{p} E^{n-1}_{X^{\mathcal{W}}}(W) + \bar F ^{p}
      E^{n-1}_{X^{\mathcal{W}}}(W))\cap 
      E^{n-1}_{X^{\mathcal{W}}}(W)_{\mathbb{R}}(p-1) \right.\\
    &\phantom{A}= E^{n-1}_{X^{\mathcal{W}},\mathbb{R}}(W)(p-1)
    \cap \bar F ^{n-p}
    E^{n-1}_{X^{\mathcal{W}}}(W))\cap 
    F^{n-p} E^{n-1}_{X^{\mathcal{W}}}(W),\\
    &\mathcal{D}_{m}(D^{X^{\mathcal{T}}}(W),q)\\
    &\phantom{A}=D_{m}^{X^{\mathcal{T}},\mathbb{R}}(W^{\infty})(q)\cap
    F_{q} D_{m}^{X^{\mathcal{T}}}(W) \cap \bar F _{q}
    D_{m}^{X^{\mathcal{T}}}(W))\\
    &\phantom{A}=D_{n-1}^{X^{\mathcal{T}},\mathbb{R}}(W)(p-1)\cap
    F_{p-1} D_{n-1}^{X^{\mathcal{T}}}(W) \cap \bar F _{p-1}
    D_{n-1}^{X^{\mathcal{T}}}(W)).
  \end{align*}
  Therefore, the first statement follows from the duality between
  $E_{X^{\mathcal{W}}}(W)$ and $D^{X^{\mathcal{T}}}(W)$ and the fact
  that, under this duality,
  $D_{n-1}^{X^{\mathcal{T}},\mathbb{R}}(W)(p-1)$ is identified
  with the dual of $E^{n-1}_{X^{\mathcal{W}},\mathbb{R}}(W)(p-1)$ and
  $F_{p-1} D_{n-1}^{X^{\mathcal{T}}}(W) $ is identified with the
  dual of $\bar F ^{n-p}
  E^{n-1}_{X^{\mathcal{T}}}(W)$.

  The compatibility with the differential is a
  straightforward computation using the formulas for the differential
  given in \cite{Burgos:CDB} theorem 2.6. For instance, 
  if $\omega \in
  \mathcal{D}^{n}(E_{X^{\mathcal{W}}}(W),p)$, with $n<2p-1$ and $T\in
  \mathcal{D}_{m}(D^{X^{\mathcal{T}}}(W),q)$, with $m=n$ and
  $q=p-1$, then we have  
  \begin{align*}
    (\dd_{\mathcal{D}} T)(\omega )
    &=(\dd T)(\omega)\\
    &=(-1)^{n}T(\dd \omega)\\
    &=(-1)^{n}T(F^{n-p+1,n-p+1}\dd \omega)\\
    &=(-1)^{n}T(-\dd_{\mathcal{D}} \omega ).
  \end{align*}
  In the third equality we have used that $T\in F_{q}\cap \bar
  F_{q}=F_{p-1,p-1}$, which implies that, for any form $\eta$, we have
  $T(\eta)=T(F^{n-p+1,n-p+1}\eta)$. The other cases are analogous.

  Similarly, the compatibility with the product follows from 
  \cite{Burgos:CDB} theorem 2.6. For instance,
  let
  $\omega \in \mathcal{D}^{n}(E_{X^{\mathcal{W}}}(W),p)$, 
  $T\in \mathcal{D}_{m}(D^{X^{\mathcal{T}}}(W),q)$ and 
  $\eta \in \mathcal{D}^{l}(E_{X^{\mathcal{W}}}(W),r)$, with
  $n-m+l=1$ and $p-q+r=1$. Assume that $n<2p$, $m>2q$, $l\ge 2r$,
  then 
  \begin{displaymath}
    (\omega \bullet T)(\eta)=
    ((-1)^{n}r_{p}(\omega )\land T+\omega \land r_{q}(T))(\eta),
  \end{displaymath}
  where $r_p(\omega )=2\pi _{p}(F^{p}\dd \omega )$ and $r_q(T)=2\pi
  _{q}(F_{q}\dd T )$. 
  But
  \begin{displaymath}
    (-1)^{n}r_{p}(\omega )\land T(\eta)=
    (-1)^{n}T(\eta \land r_{p}(\omega )),
  \end{displaymath}
  and
  \begin{align*}
    (\omega \land r_{q}(T))(\eta)&= r_{q}(T)(\eta\land \omega )\\
    &=2\pi _{q}F_{q}(\dd T)(\eta\land \omega )\\
    &= 2 F_{q}(\dd T)(\eta\land \omega )\\
    &= 2 \partial T_{q+1,m-q}(\eta\land \omega )\\
    &= T \left( 2 (-1)^{m-1}\partial (\eta\land \omega)^{q,m-q}\right)\\
    &= T \left( 2 (-1)^{n+l}\partial (\eta\land \omega)^{p+r-1,n+l-p-r}\right).
  \end{align*}
  On the other hand 
  \begin{displaymath}
    T(\eta \bullet \omega )=T\left(\eta \land r_{p}(\omega )+
    (-1)^{l}2 \partial (\omega \land \eta)^{p+r-1,n+l-p-r}\right).
  \end{displaymath}
  The other cases are analogous.
\end{proof}

\nnpar{Duality.} We summarize in the next proposition the
basic properties of formal Deligne cohomology and tempered Deligne
homology that follow from the 
previous discussions. 

\begin{proposition}\label{prop:2} For every pair of integers $n$ and
  $p$, by applying the exact functors
  $\mathcal{D}^{\ast}(\underline{\phantom{A}},p)$ and  
  $\mathcal{D}_{\ast}(\underline{\phantom{A}},p-1)$ to the
  diagrams \eqref{eq:1} and 
  \eqref{eq:2} respectively, we obtain the corresponding diagrams of
  Deligne complexes that are the topological dual of each other. In
  particular we obtain long exact sequences
  \begin{multline}\label{eq:13}
    H^{n}_{\mathcal{D}^{f},c}(W^{\infty},\mathbb{R}(p))
    \rightarrow H^{n}_{\mathcal{D}^{f}}(Y^{\infty},\mathbb{R}(p))
    \rightarrow H^{n}_{\mathcal{D}^{f}}(Z^{\infty},\mathbb{R}(p))
    \rightarrow \\ H^{n+1}_{\mathcal{D}^{f},c}(W^{\infty},\mathbb{R}(p))
    \rightarrow  
  \end{multline}
  and
  \begin{multline}\label{eq:12}
    \leftarrow H_{n-1}^{\mathcal{D}^{\mathcal{T}}}(W^{\infty},\mathbb{R}(p-1))
    \leftarrow
    H_{n-1}^{\mathcal{D}^{\mathcal{T}}}(Y^{\infty},\mathbb{R}(p-1))
    \leftarrow \\
    H_{n-1}^{\mathcal{D}^{\mathcal{T}}}(Z^{\infty},\mathbb{R}(p-1))
    \leftarrow H_{n}^{\mathcal{D}^{\mathcal{T}}}(W^{\infty},\mathbb{R}(p-1))
  \end{multline}
  and pairings 
  \begin{align*}
    H^{n}_{\mathcal{D}^{f}}(Y^{\infty},\mathbb{R}(p))\otimes
    H_{n-1}^{\mathcal{D}^{\mathcal{T}}}(Y^{\infty},\mathbb{R}(p-1))
    &\longrightarrow \mathbb{R},\\
    H^{n}_{\mathcal{D}^{f},c}(W^{\infty},\mathbb{R}(p))\otimes
    H_{n-1}^{\mathcal{D}^{\mathcal{T}}}(W^{\infty},\mathbb{R}(p-1)) 
    &\longrightarrow \mathbb{R},\\
    H^{n}_{\mathcal{D}^{f}}(Z^{\infty},\mathbb{R}(p))\otimes
    H_{n-1}^{\mathcal{D}^{\mathcal{T}}}(Z^{\infty},p-1) 
    &\longrightarrow \mathbb{R}.
  \end{align*}
  that are compatible with the above sequences.

  Moreover, the topologies of the space of differential forms and  of the
  space of currents induce structures of topological vector spaces on
  the real formal Deligne cohomology groups and the tempered Deligne
  homology groups. The 
  above pairings induce a perfect pairing of the corresponding
  separated vector spaces.
\end{proposition}
\begin{proof}
  This is a direct consequence of the exactness of the functors
  $\mathcal{D}^{\ast}(\underline{\phantom{A}},p)$ and
  $\mathcal{D}_{\ast}(\underline{\phantom{A}},p-1)$ and proposition
  \ref{prop:5}. 
\end{proof}

  The image of $\dd_{\mathcal{D}}$ in the complex
  $\mathcal{D}^{\ast}(E_{\flat}(U),p)$ does not need to be
  closed. Therefore the pairing between formal cohomology and
  tempered homology  do
  not need to be perfect. Only the induced pairing in the
  corresponding separated vector spaces is perfect. Nevertheless, in
  the case of a 
  proper algebraic complex manifold $X$, by Hodge theory, we obtain
  a perfect pairing between Deligne-Beilinson cohomology and homology.

\begin{corollary}[Exceptional duality for Deligne cohomology] \label{cor:3}
  Let $X$ be a proper complex algebraic manifold, equidimensional of
  dimension $d$. 
  Then there is a 
  perfect duality
  \begin{displaymath}
    H_{\mathcal{D}}^{n}(X,\mathbb{R}(p))\otimes
    H_{\mathcal{D}}^{2d-n+1}(X,\mathbb{R}(d-p+1))\longrightarrow \mathbb{R}
  \end{displaymath}
  which is compatible, up to a sign, with the product in Deligne
  cohomology. 
\end{corollary}
\begin{proof}
  By Poincar\'e duality in Deligne cohomology (cf. \cite{Jannsen:DcHD}
  1.5) there is a natural isomorphism
  \begin{displaymath}
    H_{\mathcal{D}}^{2d-n+1}(X,\mathbb{R}(d-p+1))\cong
    H^{\mathcal{D}}_{n-1}(X,\mathbb{R}(p-1)).
  \end{displaymath}
  By Hodge theory we know that
  \begin{displaymath}
    H_{\mathcal{D}}^{n}(X,\mathbb{R}(p))=
    \begin{cases}
      H^{n-1}(X,\mathbb{R}(p-1))\cap \overline F^{n-p} \cap F^{n-p},&
      \text{ if }  n<2p,\\
      H^{n}(X,\mathbb{R}(p))\cap \overline F^{p} \cap F^{p},&
      \text{ if }  n\ge 2p.
    \end{cases}
  \end{displaymath}
  Moreover, the pairing is given, up to a sign, by the wedge product of
  differential forms followed by the integral along $X$.
    Therefore, by Serre's duality, the pairing of
    proposition  \ref{prop:2} is perfect.
\end{proof}

\subsection{Semi-purity of tempered Deligne cohomology }
\label{sec:purity-form-deligne}
\ 

\nnpar{Vanishing theorems.} The aim of this section is to prove the
following result

\begin{theorem}(Semi-purity of tempered Deligne homology)
  Let $X$ be a projective complex algebraic manifold,
  $W$ a locally closed subvariety, of dimension at most $p$. Then
  \begin{displaymath}
    H_{n}^{\mathcal{D}^{\mathcal{T}}}(W^{\infty},\mathbb{R}(e))=0, \text{ for
    all } n > \max(e+p,2p-1). 
  \end{displaymath}
\end{theorem}
\begin{proof}
  We will prove the result by ascending induction over $p$. The
  result is trivially true for $p<0$. Then, by the exact sequence
  \eqref{eq:12} and
  induction, one is  
   reduced to the case $W$ closed. 

  We will deduce the theorem by duality from the following proposition
  \begin{proposition}\label{prop:1}
    Let $Y$ be a closed subvariety of a projective complex algebraic
    manifold. Let $p$ be the dimension of $Y$. Then
    \begin{displaymath}
      H^{n+1}_{\mathcal{D}^{f}}(Y^{\infty},\mathbb{R}(e+1))=0, \text{ for
        all } n > max(e+p,2p-1)  
    \end{displaymath}
  \end{proposition}
  \begin{proof}
    Let $\mathscr{I}_{Y}$ be the ideal of holomorphic functions on $X$
    vanishing at $Y$. We denote
    \begin{displaymath}
      \Omega ^{q}_{Y^{\infty}}= \lim_{\substack{\longleftarrow\\k}}
      \Omega ^{q}_{X}\left/ \mathscr{I}_{Y}^{k}\Omega ^{q}_{X}. \right.
    \end{displaymath}
    By \cite{KashiwaraSchapira:mfcacs} theorem 5.12 we have
    \begin{lemma} \label{lemm:1}
      The complex of sheaves $\mathscr{E}^{q,\ast}_{Y,\mathbb{R}}$ is a
      fine resolution of $\Omega ^{q}_{Y^{\infty}}$.
    \end{lemma}

    
    Since, by \cite{Poly:shcsesa}, the sheaf
    $\mathscr{E}^{\ast}_{Y^{\infty},\mathbb{R}}$ is an acyclic
    resolution of the constant sheaf $\underline{\mathbb{R}}_{Y}$,
    from lemma \ref{lemm:1} and the techniques of \cite{Burgos:CDB},
    we deduce that 
    $H^{\ast}_{\mathcal{D}^{f}}(Y^{\infty},\mathbb{R}(e+1))$ is
    isomorphic to the hypercohomology of the complex of sheaves 
    \begin{equation}\label{eq:3}
      \underline{\mathbb{R}}_{\mathcal{D}^{f},Y^{\infty}}(e):=
      \underline{\mathbb{R}}_{Y}(f)\longrightarrow \Omega  
      ^{0}_{Y^{\infty}}\longrightarrow  \dots \longrightarrow 
      \Omega ^{e}_{Y^{\infty}}.
    \end{equation}
    
    \begin{lemma} \label{lemm:2} If $n> p$ then $H^{n}(Y,\Omega _{Y^{\infty}}^{q})=0$.
    \end{lemma}
    \begin{proof}
      By \cite{hartshorne75:Rhamcag} proposition I.6.1
      \begin{displaymath}
        H^{n}(Y,\Omega _{Y^{\infty}}^{q})=
        H^{n}(Y^{\alg},\hat \Omega _{Y}^{q}),
      \end{displaymath}
      where $Y^{\alg}$ is the corresponding algebraic variety and
      $\hat \Omega _{Y}^{q}$ is the completion of the sheaf of 
      algebraic differentials. But now $Y^{\alg}$ is a noetherian
      topological space of dimension $p$, hence the lemma.
    \end{proof}

    Using lemma \ref{lemm:2} we obtain that the $E_{1}^{s,t}$ term of the
    spectral sequence of the hypercohomology of the complex
    \eqref{eq:3} can be non zero only for $s=0$, $0\le t \le 2p$ and 
    $1\le s \le e+1$, $0\le t \le p$, which implies proposition
    \ref{prop:1}. 
  \end{proof}
  
  We finish now the proof of the theorem. By proposition \ref{prop:1},
  for every $n > \max(p+e,2p-1)$, the morphism
  \begin{displaymath}
    \dd_{\mathcal{D}}^{n}: \mathcal{D}^{n}(E_{X^{\mathcal{W}}}(Y),e+1)
    \longrightarrow \mathcal{D}^{n+1}(E_{X^{\mathcal{W}}}(Y),e+1)
  \end{displaymath}
  satisfies $\Img
  (\dd_{\mathcal{D}}^{n})=\Ker(\dd_{\mathcal{D}}^{n+1})$, hence the
  image of $\dd_{\mathcal{D}}^{n}$ is
  a closed subspace. Therefore, by
  \cite{bourbaki87:_topol_vector_spaces_chapt} IV.2 theorem 1, we have
  that the dual morphism
  \begin{displaymath}
    \dd_{\mathcal{D}}:\mathcal{D}_{n}(D^{X^{\mathcal{T}}}(Y),e) \longrightarrow 
    \mathcal{D}_{n-1}(D^{X^{\mathcal{T}}}(Y),e)
  \end{displaymath}
  has closed image. This implies that, for $n\ge \max(p+e,2p-1)$, the
  vector space
  $H_{n}^{\mathcal{D}^{\mathcal{T}}}(Y^{\infty},\mathbb{R}(e))$ is separated. 
  Therefore, by proposition \ref{prop:2},
  for $n>\max(p+e,2p-1)$ the pairing 
  \begin{displaymath}
    H^{n+1}_{\mathcal{D}^{f}}(Y^{\infty},\mathbb{R}(e+1))\otimes 
    H_{n}^{\mathcal{D}^{\mathcal{T}}}(Y^{\infty},\mathbb{R}(e))
    \longrightarrow \mathbb{R}
  \end{displaymath}
  is perfect. Hence by proposition \ref{prop:1} we obtain the theorem.
\end{proof}

\nnpar{semi-purity of tempered Deligne cohomology.}
The semi-purity theorem can be stated  in terms of tempered Deligne
cohomology as follows.
\begin{corollary} \label{cor:1}
  Let $X$ be a complex quasi-projective manifold  and $Y$ a closed
  subvariety of codimension at least $p$. Then 
  \begin{displaymath}
    H^{n}_{\mathcal{D}^{\mathcal{T}}, Y}(X,\mathbb{R}(e))=0, \text{ for
      all } n < \min(e+p,2p+1), 
  \end{displaymath}
  In particular
  \begin{displaymath}
    H^{n}_{\mathcal{D}^{\mathcal{T}}, Y}(X,\mathbb{R}(p))=0, \text{ for
      all } n < 2p. 
  \end{displaymath}
\end{corollary}
This is the weak purity property used in \cite{BurgosKramerKuehn:cacg} 6.4.

\section{Arithmetic Intersection Theory}
\label{sec:arithm-inters-theory}

\subsection{ Definition of Covariant arithmetic Chow groups}
\label{sec:covar-arithm-chow}
In \cite{Burgos:acr}, the author introduced a variant of the
arithmetic Chow groups that are covariant with respect to arbitrary
proper morphisms. 
In the paper \cite{BurgosKramerKuehn:cacg} these groups are further
studied as an example of cohomological arithmetic Chow groups. These
groups are denoted by $\cha^{\ast}
(X,\mathcal{D}_{\D})$. The semi-purity property (corollary
\ref{cor:1})  was announced in \cite{BurgosKramerKuehn:cacg} and has
consequences in the behavior of the covariant arithmetic Chow
groups. On the other hand, Kawaguchi and Moriwaki
\cite{KawaguchiMoriwaki:isfav} have
given another definition of covariant arithmetic Chow
groups called $D$-arithmetic Chow groups. A consequence of Corollary
\ref{cor:1} is that, when $X$ is 
equidimensional and generically projective, both definitions of
covariant arithmetic Chow 
groups agree. We note that Zha \cite{zha99:_rieman_roch} has also
introduced a notion of covariant arithmetic Chow
groups that only differs from the definition of
\cite{KawaguchiMoriwaki:isfav} on the fact that he neglects the
anti-linear involution $F^{\infty}$.

In this section we will summarize the properties of the
covariant arithmetic Chow groups. We will follow the notations and 
terminology of
\cite{BurgosKramerKuehn:cacg}, but we will use the grading by
dimension that is more natural when dealing with covariant Chow
groups.

\nnpar{Arithmetic rings and arithmetic varieties.} Let $A$ be an
arithmetic ring (see \cite{GilletSoule:ait}) with fraction field
$F$. In particular
$A$ is provided with a non empty set of complex embeddings $\Sigma $
and a conjugate linear involution $F_{\infty}$ of $\mathbb{C}^{\Sigma
}$ that commutes with the diagonal embedding of $A$ in $\mathbb{C}^{\Sigma
}$. 
 Since we will be working with 
dimension of cycles, following \cite{GilletSoule:aRRt} we
will further impose that $A$ is equicodimensional and Jacobson.
Let $S=\Spec A$ and let $e=\dim S$.

An arithmetic variety $X$ is
a flat quasi-projective scheme over $A$, that has smooth generic fiber
$X_{F}$. To every arithmetic variety $X$ we can associate a complex
algebraic manifold $X_{\Sigma }$ and a real algebraic manifold
$X_{\mathbb{R}}=(X_{\Sigma },F_{\infty})$. 

\nnpar{The arithmetic complex of tempered Deligne homology.} To every
pair of integers $n,p$, and 
every open Zariski subset $U$ of $X_{\mathbb{R}}$ we assign the group
\begin{displaymath}
  \mathcal{D}_{n}^{\D,X}(U,p)=\mathcal{D}_{n}\left(
    D_{\ast}^{X_{\Sigma }^{\mathcal{T}}}(U),p\right)^{\sigma}, 
\end{displaymath}
where $\sigma $ is the involution that acts as complex conjugation on
the space and on the currents. That is, if $T\in
D_{n}(X_{\mathbb{C}})$ then $\sigma 
(T)=\overline {(F_{\infty})_{\ast}T}$. 
And $(\phantom{A})^{\sigma}$ denote the elements that are fixed by
$\sigma $. Then $\mathcal{D}_{n}^{\D,X}(\underline{\phantom{A}},p)$ is
a totally acyclic sheaf (in the sense of
\cite{BurgosKramerKuehn:cacg}) for the real scheme underlying
$X_{\mathbb{R}}$. When $X$ is fixed, $\mathcal{D}_{\ast}^{\D,X}$ will
be denoted by $\mathcal{D}_{\ast}^{\D}$.

If $U$ is a Zariski open subset of $X_{\mathbb{R}}$ and
$Y=X\setminus U_{\mathbb{R}}$ we write
\begin{align}
  H^{\mathcal{D}^{\mathcal{T}}}_{\ast}(U,\mathbb{R}(p))&=
  H_{\ast}(\mathcal{D}^{\D}(U,p)),\\
  H^{\mathcal{D}^{\mathcal{T}},Y}
  _{\ast}(X_{\mathbb{R}},\mathbb{R}(p))&=  
  H_{\ast}(s(\mathcal{D}^{\D}(U,p),
  \mathcal{D}^{\D}(X_{\mathbb{R}},p))),\\  
  \widetilde {\mathcal{D}}_{2p-1}^{\D}(X_{\mathbb{R}},p)&=
  \mathcal{D}_{2p-1}^{\D}(X_{\mathbb{R}},p)\left/
    \Img \dd_{\mathcal{D}},\right.\\
  {\rm Z} \mathcal{D}_{2p}^{\D}(X_{\mathbb{R}},p)&=
  \Ker(\dd_{\mathcal{D}}:\mathcal{D}_{2p}^{\D}(X_{\mathbb{R}},p)\longrightarrow 
  \mathcal{D}_{2p+1}^{\D}(X_{\mathbb{R}},p)).
\end{align}

Let $\mathcal{Z}_{p}=\mathcal{Z}_{p}(X_{\mathbb{R}})$ be the set of
dimension $p$ Zariski closed 
subsets of $X_{\mathbb{R}}$ ordered by
inclusion. Then we will write 
\begin{align*}
  \mathcal{D}_{\ast}^{\D}(X_{\mathbb{R}}\setminus \mathcal{Z}_{p},p)
  &=\lim_{\substack{\longrightarrow 
      \\ Y\in\mathcal{Z}^{p}}}
  \mathcal{D}_{\ast}^{\D}(X_{\mathbb{R}}\setminus Y,p),\\
  \widetilde{\mathcal{D}}_{\ast}^{\D}(X_{\mathbb{R}}\setminus
  \mathcal{Z}_{p},p) 
  &= \mathcal{D}_{\ast}^{\D}(X_{\mathbb{R}}\setminus
  \mathcal{Z}_{p},p)\left /
    \Img \dd_{\mathcal{D}}\right.,\\
  H^{\mathcal{D}^{\mathcal{T}},\mathcal{Z}_{p}}_{\ast}
  (X_{\mathbb{R}},\mathbb{R}(p))&=  
  H_{\ast}(s(\mathcal{D}^{\D}(X_{\mathbb{R}}\setminus
  \mathcal{Z}_{p},p), 
  \mathcal{D}^{\D}(X_{\mathbb{R}},p))).
\end{align*}

\nnpar{Green objects.} We recall the definition of Green object for a
cycle given in \cite{BurgosKramerKuehn:cacg} but adapted to the grading by
dimension.  Let $y$ be a dimension $p$ algebraic cycle of
$X_{\mathbb{R}}$. Let $Y$ be the support of $y$. The class of $y$ in
$H^{\mathcal{D}^{\mathcal{T}},Y}_{2p}(X_{\mathbb{R}},\mathbb{R}(p))$,
denoted $\cl(y)$,
is represented by the pair $(\delta _{y},0)\in
s(\mathcal{D}^{\D}(X_{\mathbb{R}},p), 
  \mathcal{D}^{\D}(U_{\mathbb{R}},p))$.
We denote also by $\cl(y)$ the image of this class in
$H^{\mathcal{D}^{\mathcal{T}},\mathcal{Z}_{p}}_{2p} 
  (X_{\mathbb{R}},\mathbb{R}(p))$.

In this setting, the truncated homology classes can be written as
\begin{multline*}
    \widehat {H}^{\mathcal{D}^{\mathcal{T}},\mathcal{Z}_{p}}_{\ast}
  (X_{\mathbb{R}},\mathbb{R}(p))=\\ \{(\omega _{y},\widetilde g_{y})\in {\rm Z}
\mathcal{D}_{2p}^{\D}(X,p)\oplus
\widetilde{\mathcal{D}}_{2p-1}^{\D}(X_{\mathbb{R}}\setminus 
  \mathcal{Z}_{p},p)\mid \dd_{\mathcal{D}} \widetilde g_{y}=\omega
  _{y}\}.
\end{multline*}

There is an obvious class map
\begin{displaymath}
  \cl: \widehat {H}^{\mathcal{D}^{\mathcal{T}},\mathcal{Z}_{p}}_{\ast} 
  (X_{\mathbb{R}},\mathbb{R}(p))\longrightarrow 
  H^{\mathcal{D}^{\mathcal{T}},\mathcal{Z}_{p}}_{\ast}
  (X_{\mathbb{R}},\mathbb{R}(p)).
\end{displaymath}
Then a Green object for $y$ is an element 
$$\mathfrak{g}_{y}=(\omega _{y},\widetilde g_{y})\in 
\widehat {H}^{\mathcal{D}^{\mathcal{T}},\mathcal{Z}_{p}}_{2p}
  (X_{\mathbb{R}},\mathbb{R}(p))
$$
such that $\cl(\mathfrak{g}_{y})=\cl(y)$.

The following result follows directly from the definition
\begin{lemma} \label{lemm:3}
  An element $\mathfrak{g}_{y}=(\omega _{y},\widetilde g_{y})\in
  \widehat {H}^{\mathcal{D}^{\mathcal{T}},\mathcal{Z}_{p}}_{2p}
  (X_{\mathbb{R}},\mathbb{R}(p))$ is a Green object for $y$ if
  and only if there exists a current $\widetilde \gamma \in \widetilde
  {\mathcal{D}}_{2p-1}^{\D}(X_{\mathbb{R}},p)$ such that 
  \begin{align*}
    \widetilde g_{y}&=\widetilde \gamma|_{X\setminus
      \mathcal{Z}_{p}}\\ 
    \dd_{\mathcal{D}} \widetilde \gamma +\delta _{y}&=\omega _{y}.
  \end{align*}
\end{lemma}

\nnpar{Arithmetic Chow groups.} Every dimension $p$ algebraic cycle
$y$ on $X$ defines a dimension $(p-e)$ algebraic cycle $y_{\mathbb{R}}$
on $X_{\mathbb{R}}$, where $e$ is the dimension of the base scheme $S$. 

\begin{definition}
  The group of arithmetic cycles of dimension $p$ is defined as
  \begin{displaymath}
    \za_{p}(X,\mathcal{D}^{\D})=
    \{(y,\mathfrak{g}_{y})\in {\rm Z}_{p}(X)\oplus
    \widehat {H}^{\mathcal{D}^{\mathcal{T}},\mathcal{Z}_{p-e}}_{2p-2e}
  (X_{\mathbb{R}},\mathbb{R}(p-e))\mid
  \cl(y_{\mathbb{R}})=\cl(\mathfrak{g}_{y})\}. 
  \end{displaymath}
  Let $W$ be a dimension $p+1$ irreducible subvariety of $X$ and $f\in
  K(W)^{\ast}$ be a rational function. Let $\widetilde W_{\mathbb{R}}$
  be a resolution of singularities of $W_{\mathbb{R}}$ and let
  $\iota:\widetilde W_{\mathbb{R}}\longrightarrow X_{\mathbb{R}}$ be
  the induced map. Then we write
  \begin{displaymath}
    \diva f = (\dv f, (0,\iota_{\ast}(-\frac{1}{2} \log f\bar f )).
  \end{displaymath}
  The group of cycles rationally equivalent to zero is the subgroup 
  $$\rata_{p}(X,\mathcal{D}^{\D})\subset \za_{p}(X,\mathcal{D}^{\D})$$
  generated by the elements 
  of the form $\diva f$. 
  The \emph{homological arithmetic Chow groups} of $X$ are defined as
  \begin{displaymath}
    \cha_{p}(X,\mathcal{D}^{\D})=\za_{p}(X,\mathcal{D}^{\D})\left /
      \rata_{p}(X,\mathcal{D}^{\D})\right.
  \end{displaymath}
\end{definition}

There are well-defined maps
\begin{alignat*}{2}
\zeta&:\cha_{p}(X,\mathcal{D}^{\D})\longrightarrow\CH_{p}(X),&&\quad 
\zeta[y,\mathfrak{g}_{y}]=[y], \\  
\rho&:\CH_{p,p+1}(X)\longrightarrow
H_{2p-2e+1}^{\mathcal{D}^{\mathcal{T}}}(X,p-e) 
\subseteq\widetilde{\mathcal{D}}_{2p+1}^{\D}(X,p),&&\quad\rho[f]= 
\cl(f), \\
\amap&:\widetilde{\mathcal{D}}_{2p-2e+1}(X,p-e)\longrightarrow\cha_{p}
(X,\mathcal{D}^{\D}),&&\quad\amap(\widetilde{a})=[0,\amap(\widetilde{a})], \\
\omega&:\cha_{p}(X,\mathcal{D}^{\D})\longrightarrow{\rm
  Z}\mathcal{D}_{2p-2e}^{\D}(X,p-e),  
&&\quad\omega[y,\mathfrak{g}_{y}]=\omega(\mathfrak{g}_{y}), \\
h&:{\rm Z}\mathcal{D}_{2p}^{\D}(X,p)\longrightarrow
H_{2p}^{\mathcal{D}^{\mathcal{T}}}(X,p),  
&&\quad h(\alpha)=[\alpha]. 
\end{alignat*}  

\subsection{ Properties of Covariant arithmetic Chow groups}
\label{sec:prop-covar-arithm}
\ 

\nnpar{Basic properties.} Recall that in
\cite{BurgosKramerKuehn:cacg}, there are defined
contravariant arithmetic Chow groups denoted by $\cha^
{\ast}(X,\mathcal{D}_{\log})$.
The following result follows from the theory developed
\cite{BurgosKramerKuehn:cacg} and corollary \ref{cor:1} (semi-purity
property).  

\begin{theorem}
\label{thm:logD}
With the above notations, we have the following statements:
\begin{enumerate}
\item[(i)] There are exact sequences
\begin{displaymath}
\CH_{p,p+1}(X)\overset{\rho}{\longrightarrow}\widetilde 
{\mathcal{D}}_{2p-2e+1}^{\D}(X,p-e)\overset{\amap}{\longrightarrow} 
\cha_{p}(X,\mathcal{D}^{\D})\overset{\zeta}{\longrightarrow}
\CH_{p}(X)\longrightarrow 0. 
\end{displaymath}
\begin{align*}
&\CH_{p,p+1}(X)\overset{\rho}{\longrightarrow}H_{2p-2e+1}^{\mathcal
{D}^{\mathcal{T}}}(X_{\mathbb{R}},\mathbb{R}(p-e))\overset{\amap}{\longrightarrow}
\cha_{p}(X,\mathcal{D}^{\D})\overset{(\zeta,-\omega)}
{\longrightarrow} \\ 
&\phantom{CH_{p,p+1}}\CH_{p}(X) 
\oplus{\rm Z}\mathcal{D}_{2p-2e}^{\D}(X,p-e)\overset{\cl+h}
{\longrightarrow}H_{2p-2e}^{\mathcal{D}^{f}}(X_{\mathbb{R}},\mathbb{R}(p-e)) 
\longrightarrow 0.
\end{align*}
In particular, if $X_{F}$ is projective, then there is an
exact sequence 
\begin{align*}
&\CH_{p,p+1}(X)\overset{\rho}{\longrightarrow}H_{2p-2e+1}^{\mathcal
{D}}(X_{\mathbb{R}},\mathbb{R}(p-e))\overset{\amap}{\longrightarrow}
\cha_{p}(X,\mathcal{D}^{\D})\overset{(\zeta,-\omega)}
{\longrightarrow} \\ 
&\phantom{CH_{p,p+1}}\CH_{p}(X) 
\oplus{\rm Z}\mathcal{D}_{2p-2e}^{\D}(X,p-e)\overset{\cl+h}
{\longrightarrow}H_{2p-2e}^{\mathcal{D}}(X_{\mathbb{R}},\mathbb{R}(p-e)) 
\longrightarrow 0.
\end{align*}
\item[(ii)]
For any regular arithmetic variety $X$ over $A$ there are defined
contravariant arithmetic Chow groups
$\cha^{p}(X,\mathcal{D}_{\log})$. Furthermore, if $X$ is
equidimensional of dimension $d$, then  there is 
a morphism of arithmetic Chow groups
\begin{displaymath}
\cha^{p}(X,\mathcal{D}_{\log})\longrightarrow \cha_{d-p}
(X,\mathcal{D}^{\D}).
\end{displaymath}
When $X_{F}$ is projective this morphism is a monomorphism. 
Moreover, if $X_{F}$ has dimension zero, this morphism is an
isomorphism.  
\item[(iii)] \label{item:1}
For any proper morphism $f:X\longrightarrow Y$ of arithmetic 
varieties over $A$, there is a 
 morphism of covariant arithmetic Chow groups
\begin{displaymath}
f_{\ast}:\cha_{p}(X,\mathcal{D}_{\D})\longrightarrow\cha_{p}
(Y,\mathcal{D}_{\D}).    
\end{displaymath}
If $g:Y\longrightarrow Z$ is another such morphism, the equality
$(g\circ f)_{\ast}=g_{\ast}\circ f_{\ast}$ holds. Moreover, if $X$ and
$Y$ are regular and 
$f_{F}:X_{F}\longrightarrow Y_{F}$ is
a smooth proper morphism of projective varieties, then $f_{\ast}$ is 
compatible with the direct image of contravariant arithmetic Chow 
groups.
\item [(iv)] If $f:X\longrightarrow Y$ is a flat morphism,
  equidimensional of relative
  dimension $d$, and such that $f_{F}$ is smooth, then there
  is a pull-back map 
  \begin{displaymath}
    f^{\ast}:\cha_{p}(Y,\mathcal{D}^{\D})\longrightarrow
    \cha_{p+d}(X,\mathcal{D}^{\D}). 
  \end{displaymath}
  If $X$ and $Y$ are regular and equidimensional, this map is equivalent with
  the pullback map defined in the contravariant Chow groups.
  \item [(v)] Let $f:X\longrightarrow Y$ be a flat map between arithmetic
  varieties, which is smooth over $F$ and let $g:P\longrightarrow Y$
  be a proper map. Let $Z$ be the fiber product of $X$ and $P$ over
  $Y$, with $p:Z\longrightarrow P$ and $q:Z\longrightarrow X$ the two
  projections. Thus $p$ is flat and smooth over $F$ and $q$ is
  proper. Then for any $x\in \cha_{\ast}(P,\mathcal{D}^{\D})$, it holds 
  \begin{displaymath}
    q_{\ast}p^{\ast}(x)=f^{\ast}g_{\ast}(x)\in
    \cha_{\ast}(X,\mathcal{D}^{\D}).
  \end{displaymath}
\end{enumerate}
\end{theorem}
\begin{proof}
  Part (i) follows from the standard exact sequences of
  \cite{BurgosKramerKuehn:cacg} Theorem 4.13 adapted to the
  grading by dimension and corollary \ref{cor:1}.

  For (ii) we first note that, if $M$ is an equidimensional complex
  algebraic 
  manifold, $D\subset X$ is a normal crossing divisor, $\omega$ is a
  differential form with logarithmic singularities along $D$ and
  $\eta$ is a form that is flat along $D$, then $\eta\wedge \omega $ is
  flat along $D$. In particular, if $M$ is proper and $U=M\setminus D$,
  then the associated 
  current $[\omega]$ belongs to $ D^{\extd}_{\ast}(U)$. Therefore, if $y$ is a
  codimension $p$ cycle on $X$ then, by the assumptions on $X$ and on the
  arithmetic ring, $y$ is a dimension $d-p$ algebraic algebraic
  cycle. Moreover, if $(\omega_{y},\widetilde g_{y})$ is a Green form
  for $y$ (i.e. a $\mathcal{D}_{\log}$-Green object for $y$) then, by
  lemma \ref{lemm:3} and
  \cite{BurgosKramerKuehn:cacg} Proposition 6.5 we have that $([\omega
  _{y}],[\widetilde g_{y}])$ is a $\mathcal{D}^{\D}$-Green object for
  $y$. Thus we have a well defined map
  \begin{displaymath}
    \za^{p}(X,\mathcal{D}_{\log})\longrightarrow
    \za_{d-p}(X,\mathcal{D}^{\D}). 
  \end{displaymath}
  By definition this map is compatible with rational equivalence,
  hence we obtain a map at the level of Chow groups.

  To prove (iii) we first observe that, if $Z\subset X_{\Sigma }$ is a
  closed subset, then $f_{\ast} D_{\ast}^{X_{\Sigma }^{\mathcal{T}}}(Z)\subset
  D^{X_{\Sigma }^{\mathcal{T}}}(f(Z))$. Therefore, the push-forward of
  currents define a 
  covariant $f$-morphism
  \begin{displaymath}
    f_{\#}:f_{\ast}\mathcal{D}_{\ast}^{\D,X} \longrightarrow
    \mathcal{D}^{\D,Y}_{\ast}. 
  \end{displaymath}
  Here we are using the terminology of \cite{BurgosKramerKuehn:cacg}
  3.67 but adapted to the grading by dimension. Therefore applying
  \cite{BurgosKramerKuehn:cacg} \S 4.5 we obtain the push-forward map
  for covariant arithmetic Chow groups.
  More concretely this map is defined as  
  \begin{displaymath}
    f_{\ast}(y,(\omega _{y},\widetilde g_{y}))=
    (f_{\ast}y,(f_{\ast} \omega
    _{y},(f_{\ast}g_{y})\widetilde{\phantom A})).
  \end{displaymath}
  It is straightforward to check that it is compatible with the
  direct image of $\mathcal{D}_{\log}$-arithmetic Chow groups when $Y$
  is projective and $f_{F}$ smooth.

  We now prove (iv). Since $f_{F}$ is smooth, for any Zariski closed subset
  $Z\subset Y_{\mathbb{R}}$ equidimensional of dimension $p$, there is a
  well defined morphism 
  $f^{\ast}D_{n}(Y_{\Sigma })\longrightarrow D_{n+2d}(X_{\Sigma })$
  that sends $D_{n}^{Y_{\Sigma }^{\mathcal{T}}}(Z)$ to $D^{X_{\Sigma
    }^{\mathcal{T}}}_{n+2d}(f^{-1}(Z))$. Therefore we obtain well
  defined morphisms  
  \begin{displaymath}
    \begin{matrix}
      f^{\#}:\mathcal{D}^{\D}_{n}(Y_{\mathbb{R}},p)&\longrightarrow&
      \mathcal{D}^{\D}_{n+2d}(X_{\mathbb{R}},p+d),\\
      f^{\#}:\mathcal{D}^{\D}_{n}(Y_{\mathbb{R}}\setminus Z,p)&\longrightarrow&
      \mathcal{D}^{\D}_{n+2d}(X_{\mathbb{R}}\setminus f^{-1}Z,p+d),
    \end{matrix}
  \end{displaymath}
  that send $T$ to $f^{\ast}T/(2\pi i)^{d}$.
  Then the proof of (iv) is straightforward using the theory of
  \cite{BurgosKramerKuehn:cacg} 4.4 adapted to the grading by
  dimension.   

  (v) Follows as \cite{GilletSoule:aRRt} Lemma 11.
\end{proof}

\nnpar{Multiplicative properties.} In the next result we state the
multiplicative properties between covariant and contravariant Chow
groups. The proofs are simple modification of \cite{GilletSoule:aRRt} 
Theorem 3. First, for a form $\eta \in
\widetilde{\mathcal{D}}^{2p-1}_{\log}(X_{\mathbb{R}},p)$ and an element 
$x\in \cha_{q}(X,\mathcal{D}^{\D})$ we define
\begin{displaymath}
  \eta \cap x = \amap(\eta \bullet \omega (x))=\amap(\eta \wedge
  \omega (x)). 
\end{displaymath}

\begin{theorem} \label{thm:1}
  Given a map $f:X\longrightarrow Y$ of arithmetic varieties,
  with $Y$ regular, there is a cap product
  \begin{displaymath}
    \begin{matrix}
      \cha^{p}(Y,\mathcal{D}_{\log})\otimes \cha_{q}(X,\mathcal{D}^{\D})
      &\longrightarrow &\cha_{q-p}(X,\mathcal{D}^{\D})_{\mathbb{Q}}\\
      y\otimes x &\longmapsto & y._{f}x
    \end{matrix}
  \end{displaymath}
  which is also denoted $y\cap X$ if $X=Y$. This product satisfies the
  following properties
  \begin{enumerate}
  \item $\omega (y._{f}x)=f^{\ast}\omega (y)\land \omega (x)$, and,
    for any $\eta \in \widetilde
    {\mathcal{D}}_{\log}^{2p-1}(Y_{\mathbb{R}},p)$, it holds 
    $\amap(\eta)._{f}x=\amap(f^{\ast}(\eta))\cap x$.
  \item $\cha_{\ast}(X,\mathcal{D}^{\D})_{\mathbb{Q}}$ is a graded
    $\cha^{\ast}(Y,\mathcal{D}_{\log})$-module.
  \item If $g:Y\longrightarrow Y'$ is a map of arithmetic varieties
    with $Y'$ also regular, $y'\in \cha^{p}(Y',\mathcal{D}_{\log})$ and
    $x\in \cha_{q}(X,\mathcal{D}^{\D})$, then
    $y'._{gf}x=g^{\ast}(y')._{f}x$.
  \item If $h:X'\longrightarrow X$ is a projective morphism, $x'\in
    \cha_{q}(X',\mathcal{D}^{\D})$ and $y\in
    \cha^{p}(Y,\mathcal{D}_{\log})$, then
    $y._{f}(h_{\ast}(x'))=h_{\ast}(y._{fh} x')$.
  \item If $h:X'\longrightarrow X$ is flat and smooth over $F$, $x\in
    \cha_{q}(X,\mathcal{D}^{\D})$, $y\in
    \cha^{p}(Y,\mathcal{D}_{\log})$, then
    $h^{\ast}(y._{f}x)=y._{f}(h^{\ast}(x))$.
  \item Let $f:X\longrightarrow Y$ be a flat map between arithmetic
  varieties, with $Y$ regular and projective, and let $g:P\longrightarrow Y$
  be a proper smooth map of arithmetic varieties of relative dimension
  $d$. Let $Z$ be the fiber product of $X$ and $P$ over
  $Y$, with $p:Z\longrightarrow P$ and $q:Z\longrightarrow X$ the two
  projections. Then, for all $x\in \cha_{p}(X,\mathcal{D}^{\D})$ and 
  $\gamma \in \cha^{q}(P,\mathcal{D}_{\log})$, it holds the equality
  \begin{displaymath}
    q_{\ast}(\gamma ._{p}q^{\ast}(x))=g_{\ast}\gamma ._{f} \alpha .
  \end{displaymath}
  \end{enumerate}
\end{theorem}
\begin{proof}
  To define $y._{f}x$ we follow closely \cite{GilletSoule:aRRt}. We
  may assume that $Y$ is equidimensional, that
  $x=(V,\mathfrak{g}_{V})$ with $V$ a prime algebraic cycle and
  $y=(W,\mathfrak{g}_{W})$ with each component of $W$ meeting $V$
  properly on the generic fiber $X_{F}$. As in \cite{GilletSoule:aRRt}
  we can define a cycle $[V]._{f}[W]\in CH_{q-p}(V\cap
  f^{-1}(|W|))_{\mathbb{Q}}$ that gives us a well defined cycle
  $([V]._{f}[W])_{F}\in {\rm Z}_{q-p}(X_{F})$. Our task now is to
  construct the Green object for this cycle. Let
  $\mathfrak{g}_{W}=(\omega _{W},\widetilde g_{W})$ and
  $\mathfrak{g}_{V}=(\omega _{V},\widetilde g_{V})$. We write
  $U_{V}=X_{\mathbb{R}\setminus |V|}$,
  $U_{W}=X_{\mathbb{R}\setminus f^{-1}|W|}$ and $r=q-p$. 

We now define,
  in analogy with \cite{BurgosKramerKuehn:cacg} theorem 3.37, 
  \begin{align*}
    \mathfrak{g}_{W}&\ast_{f} \mathfrak{g}_{V}=
    f^{\ast}\mathfrak{g}_{W}\ast \mathfrak{g}_{V}\\
    &=\left(f^{\ast}(\omega_{W})\bullet\omega_{V},
      ((f^{\ast}(g_{W})\bullet\omega_{V},f^{\ast}(\omega_{W})\bullet
      g_{V}),-f^{\ast}(g_{W})\bullet g_ 
      {V})^{\widetilde{\phantom{=}}}\right)\\
    &=(f^{\ast}(\omega_{W})\wedge \omega_{V},
      ((f^{\ast}(g_{W})\wedge \omega_{V},f^{\ast}(\omega_{W})\wedge
      g_{V}),\\
    & \ 
      \partial f^{\ast}(g_{W})\land
      g_{V}-\bar{\partial}f^{\ast}(g_{W})\land g_{V}- 
      f^{\ast}(g_{W})\land\partial g_{V}+f^{\ast}(g_{W})\land\bar{\partial}g_{V}
      )^{\widetilde{\phantom{=}}}\\
    &\in \widehat H_{2e}(\mathcal{D}^{\D}_{\ast}(X_{\mathbb{R}},e),
    s(\mathcal{D}^{\D}_{\ast}(U_{W},e)\oplus
    \mathcal{D}^{\D}_{\ast}(U_{V},e)\rightarrow
    \mathcal{D}^{\D}_{\ast}(U_{W}\cap U_{V},e)))\\
    &\cong  \widehat H_{2e}(\mathcal{D}^{\D}_{\ast}(X_{\mathbb{R}},e),
    \mathcal{D}^{\D}_{\ast}(U_{W}\cup U_{V},e)).
  \end{align*}
  Now the proof follows as in \cite{GilletSoule:aRRt} Theorem 3 and
  Lemma 12. 
\end{proof}

\begin{remark}
  \begin{enumerate}
  \item The main difference between the arithmetic Chow groups
    introduced here and the arithmetic Chow groups used in
    \cite{GilletSoule:aRRt} is that, if $x\in
    \cha_{\ast}(X,\mathcal{D}^{\D})$ then $\omega (x)$ is an arbitrary
    current instead of a smooth differential form. This allows us to
    define direct images for arbitrary proper morphisms. But the price
    we have to pay is that there ire defined inverse images only for
    morphisms that are smooth over $F$.
  \item The fact that
    the compatibility of direct images for the covariant Chow groups
    and direct images for the contravariant Chow groups in theorem \ref{thm:logD}
    \label{item:2} is stated only for varieties that are generically
    projective, is due to the fact that the
    latter is only defined when the base is proper. There are two ways
    to overcome this difficulty. One is to allow arbitrary
    singularities at infinity in the spirit of
    \cite{BurgosKramerKuehn:accavb} 3.5, but then, one will have to
    allow also arbitrary singularities at infinity for currents. This
    means that we will have to consider currents that are tempered
    in some components of the boundary but are not tempered in the
    other. The second option would be to use a different notion of
    logarithmic singularities that has better properties with respect
    to direct images. 
  \end{enumerate}

\end{remark}

\nnpar{Relationship with other arithmetic Chow groups.} Let us assume
now that $X_{F}$ is projective and let $\cha^{\ast}(X)$ denote the 
arithmetic Chow groups introduced in \cite{GilletSoule:ait} and 
$\cha_{\ast}(X)$ denote the arithmetic Chow groups introduced in
\cite{GilletSoule:aRRt}. In \cite{BurgosKramerKuehn:cacg} it is shown
that there is an isomorphism
\begin{displaymath}
  \psi :\cha^{\ast}(X,\mathcal{D}_{\log})\longrightarrow \cha^{\ast}(X),
\end{displaymath}
that is compatible with products, inverse images with respect to
arbitrary morphisms and direct images with respect to proper morphism
that are smooth over $F$. We shall state the analogous result for
covariant arithmetic Chow groups.

\begin{proposition}\label{prop:6}
  Let $X$ be an arithmetic variety with $X_{F}$ projective. Then there
  is a short exact sequence
  \begin{multline*}
    0\longrightarrow  \cha_{\ast}(X)\overset{\phi }{\longrightarrow }
    \cha_{\ast}(X,\mathcal{D}^{\D})\\
    \longrightarrow 
    \bigoplus_{p}{\rm Z}\mathcal{D}_{2p}^{\D}(X_{\mathbb{R}},p)\left/
      {\rm Z} \mathcal{D}_{2p}^{\smooth}(X_{\mathbb{R}},p)\longrightarrow
      0\right. ,
  \end{multline*}
  where $\mathcal{D}_{2p}^{\smooth}(X_{\mathbb{R}},p)$ denotes the
  subspace of currents that can be represented by smooth differential
  forms.  Moreover $\phi $ satisfies the following properties
  \begin{enumerate}
  \item If $f:X\longrightarrow Y$ is a proper morphism of arithmetic
    varieties that is smooth over $F$ and with $Y_{F}$ projective,
    then $f_{\ast}\circ \phi =\phi 
    \circ f_{\ast}$.
  \item If $f:X\longrightarrow Y$ is a flat morphism of arithmetic
    varieties that is smooth over $F$, with $X_{F}$ and $Y_{F}$
    projective, then $f^{\ast}\circ \phi =\phi 
    \circ f^{\ast}$.
  \item If $f:X\longrightarrow Y$ is a morphism of arithmetic
    varieties, with $X_{F}$ and $Y_{F}$ projective and $Y$ regular
    then, for $y\in \cha^{p}(Y,\mathcal{D}_{\log})$ and $x\in
    \cha_{q}(Y)$, it holds the equality
    \begin{displaymath}
      y._{f} \phi (x) = \psi (y)._{f} x.
    \end{displaymath}
  \end{enumerate}
\end{proposition}
\begin{proof}
  Let $y$ be a dimension $p$ algebraic cycle of
$X$ and let $g_{y}$ be a Green current for $y$ in the sense of
\cite{GilletSoule:aRRt}. Recall that the
normalization used here for the current $\delta _{y}$ differs with
the normalization used in \cite{GilletSoule:aRRt} by a factor
$\frac{1}{(2\pi i)^{p}} $. Then, by \ref{lemm:3}, the pair
\begin{displaymath}
  \left(\frac{1}{2(2\pi  i)^{p+1}}g_{y}|_{X_{\mathbb{R}}\setminus
    \mathcal{Z}_{p}},
  \frac{1}{2(2\pi  i)^{p+1}} (-2\partial\bar \partial) g_{y}+\delta
  _{y} \right) 
\end{displaymath}
is a $\mathcal{D}^{\D}$-Green object for $y$. Therefore we obtain a
well defined morphism 
$\za_{p}(X)\longrightarrow \za_{p}(X,\mathcal{D}^{\D})$. It is
straightforward to check that
this map preserves rational equivalence, the exactness of the above
exact sequence and properties (i), (ii) and (iii). 
\end{proof}

\begin{corollary} With the hypothesis of the proposition, 
  every element $x\in \cha_{p}(X,\mathcal{D}^{\D})$ can be represented as
  \begin{displaymath}
    x= \phi(x_{1})+\amap(\eta) 
  \end{displaymath}
  where $x_{1}\in \cha_{p}(X)$ and $\eta\in \widetilde
  {\mathcal{D}}_{2p+1}^{\D}(X_{\mathbb{R}},p)$. Moreover, if 
  \begin{displaymath}
     x= \phi(x_{1})+\amap(\eta) = \phi(x'_{1})+\amap(\eta')
  \end{displaymath}
  are two such representations, then $\eta-\eta'\in \widetilde
  {\mathcal{D}}_{2p+1}^{\smooth}(X_{\mathbb{R}},p).$
\end{corollary}
\begin{proof}
  This follows from the previous proposition and the fact that the map
  \begin{displaymath}
    \dd_{\mathcal{D}}:\widetilde
  {\mathcal{D}}_{2p+1}^{\D}(X_{\mathbb{R}},p)\longrightarrow 
  {\rm Z}\mathcal{D}_{2p}^{\D}(X_{\mathbb{R}},p)\left/
      {\rm Z} \mathcal{D}_{2p}^{\smooth}(X_{\mathbb{R}},p)
      \right.
  \end{displaymath}
  is surjective due to the projectivity of $X$. The last statement
  follows from \cite{GilletSoule:ait} Theorem 1.2.2.
\end{proof}

The following result follows now easily from the previous corollary.

\begin{corollary}
  Assume furthermore that $X$ is equidimensional of dimension $d$ and
  let $\cha^{\ast}_{D}(X)$ denote the $D$-arithmetic Chow groups
  introduced in \cite{KawaguchiMoriwaki:isfav}. Then there is a
  natural isomorphism
  \begin{displaymath}
    \bigoplus _{p}\cha^{p}_{D}(X)\longrightarrow \bigoplus
    _{p}\cha_{d-p}(X,\mathcal{D}^{\D}). 
  \end{displaymath}
  Moreover this isomorphism is compatible with push-forwards and
  the structure of module over the contravariant arithmetic Chow
  groups. \hfill $\square$
\end{corollary}

\bibliographystyle{amsplain}

\begin{thebibliography}{10}

\bibitem{bourbaki87:_topol_vector_spaces_chapt}
N.~Bourbaki, \emph{Topological vector spaces, chapters 1-5}, Springer Verlag,
  1987.

\bibitem{BrasseletPflaum:_whitn}
J.P. Brasselet and M.~Pflaum, \emph{On the homology of algebras of {W}hitney
  functions on subanalytic sets}, Institut de Mathématiques de Luminy,
  Prétirage no 2002-02, March 2002.

\bibitem{Burgos:acr}
J.I. Burgos~Gil, \emph{Arithmetic {C}how rings}, Ph.D. thesis, University of
  Barcelona, 1994.

\bibitem{Burgos:CDB}
\bysame, \emph{Arithmetic {C}how rings and {D}eligne-{B}eilinson cohomology},
  J. Alg. Geom. \textbf{6} (1997), 335--377.

\bibitem{BurgosKramerKuehn:accavb}
J.I. Burgos~Gil, J.~Kramer, and U.~K\"uhn, \emph{Arithmetic characteristic
  classes of automorphic vector bundles}, Documenta Math. \textbf{10} (2005),
  619--716.

\bibitem{BurgosKramerKuehn:cacg}
\bysame, \emph{Cohomological arithmetic
  {C}how rings}, J. Inst. Math. Jussieu \textbf{6} (2007), no.~1, 1--172.
  \MR{MR2285241}


\bibitem{GilletSoule:ait}
H.~Gillet and C.~Soul{\'e}, \emph{Arithmetic intersection theory}, Publ. Math.
  IHES \textbf{72} (1990), 94--174.

\bibitem{GilletSoule:aRRt}
\bysame, \emph{An arithmetic {R}iemann-{R}och theorem}, Invent. Math.
  \textbf{110} (1992), 473--543.

\bibitem{GriffithsHarris:pag}
P.~Griffiths and J.~Harris, \emph{Principles of algebraic geometry}, John Wiley
  {\&} Sons, Inc., 1994.

\bibitem{hartshorne75:Rhamcag}
R.~Hartshorne, \emph{On the de {R}ham cohomology of algebraic varieties}, Publ.
  math. IHES \textbf{45} (1975), 5--99.

\bibitem{Jannsen:DcHD}
U.~Jannsen, \emph{Deligne homology, {H}odge-{D}-conjecture and motives},
  Beilinson's Conjectures on Special Values of $L$-Functions (M.~Rapoport,
  N.~Schappacher, and P.~Schneider, eds.), Perspectives in Math., vol.~4,
  Academic Press, 1988, pp.~305--372.

\bibitem{KashiwaraSchapira:mfcacs}
Masaki Kashiwara and Pierre Schapira, \emph{Moderate and formal cohomology
  associated with constructible sheaves}, M\'em. Soc. Math. France (N.S.)
  (1996), no.~64, iv+76. \MR{MR1421293 (97m:32052)}

\bibitem{KawaguchiMoriwaki:isfav}
S.~Kawaguchi and A.~Moriwaki, \emph{Inequalities for semistable families of
  arithmetic varieties}, Preprint alg-geom 9710007, 1998.

\bibitem{Poly:shcsesa}
J.P. Poly, \emph{Sur l'homologie des courants \`a support dans un ensemble
  semi- analytique}, M\'em. Soc. Math. France \textbf{38} (1974), 35--43.

\bibitem{Tougeron:Ifd}
J.C. Tougeron, \emph{Id\'eaux de fonctions diff\'erentiables}, Springer-Verlag,
  1972.

\bibitem{zha99:_rieman_roch}
Yuhan Zha, \emph{A general arithmetic {R}iemann-{R}och theorem}, Ph.D. thesis,
  University of Chicago, 1998.

\end{thebibliography}
\newcommand{\noopsort}[1]{} \newcommand{\printfirst}[2]{#1}
  \newcommand{\singleletter}[1]{#1} \newcommand{\switchargs}[2]{#2#1}
\providecommand{\bysame}{\leavevmode\hbox to3em{\hrulefill}\thinspace}
\providecommand{\MR}{\relax\ifhmode\unskip\space\fi MR }
\providecommand{\MRhref}[2]{%
  \href{http://www.ams.org/mathscinet-getitem?mr=#1}{#2}
}
\providecommand{\href}[2]{#2}

\end{document}